\documentclass[11pt]{article}
\usepackage{graphicx, caption, subcaption}
\usepackage{cite}
\usepackage{amsmath, amssymb, hyperref, multirow, longtable, paracol, booktabs, diagbox, array, makecell}
\usepackage{amsthm}
\usepackage{dsfont}
\usepackage[papersize={8.5in,11in}]{geometry}
\usepackage[normalem]{ulem}
\usepackage{xcolor}
\usepackage{mathtools}
\usepackage{mathrsfs}
\usepackage[toc,page]{appendix}
\usepackage{authblk}

\DeclarePairedDelimiter{\abs}{\lvert}{\rvert}
\DeclarePairedDelimiter{\norm}{\lVert}{\rVert}

\newtheorem{theorem}{Theorem}[section]
\newtheorem{assumption}[theorem]{Assumption}
\newtheorem{lemma}[theorem]{Lemma}
\newtheorem{proposition}[theorem]{Proposition}
\newtheorem{corollary}[theorem]{Corollary}
\newtheorem{definition}[theorem]{Definition}
\newtheorem{remark}[theorem]{Remark}

\newcommand{\R}{\mathbb{R}}
\newcommand{\I}{[0,1]}
\newcommand{\Tw}{T_{W}}
\newcommand{\inner}[2]{\left\langle #1, #2\right\rangle}
\newcommand{\dd}{\mathrm{d}}
\newcommand{\Mn}{\mathcal{M}_{n}}
\newcommand{\Pn}{P_{n}}
\newcommand{\rank}{\operatorname{rank}}
\newcommand{\ran}{\operatorname{Ran}}

\newcommand{\spec}{\operatorname{spec}}
\newcommand{\Mr}{\mathcal{M}_{r}}
\newcommand{\src}[2][]{\kern-\fontdimen2\font\cite[#1]{#2}}

\definecolor{forestgreen}{rgb}{0.33,0.61,0.34}

\title{\vspace{-1.5cm}Spectral dimension reduction of dynamical systems on graphons}



\author[1,2,*]{Naoki Masuda}
\author[2]{Bisna Mary Eldo}
\affil[1]{Gilbert S. Omenn Department of Computational Medicine and Bioinformatics, University of Michigan, Ann Arbor, MI 48109-2218, USA}
\affil[2]{Department of Mathematics, University of Michigan, Ann Arbor, MI 48109-1043, USA}
\affil[*]{Corresponding author: \href{mailto:naokimas@umich.edu}{naokimas@umich.edu}}
\date{}

\begin{document}
\maketitle

\begin{abstract}

Spectral dimension reduction approximates a dynamical system on a large network by a low-dimensional system whose variables are weighted averages of the node states, with weights given by dominant eigenvectors of the network. Its effectiveness is closely related to the low-rank structure of the interaction matrix. We develop a graphon analogue of this reduction for dynamical systems on dense graph limits. For graphon dynamical systems, we formulate the reduction in terms of the compact self-adjoint integral operator associated with the graphon. This operator-theoretic setting gives a continuum analogue of the spectral/SVD-based framework used for finite networks: the eigenfunctions provide the reduction basis, and the decay of the corresponding singular values quantifies the effect of discarded modes. For finite-rank graphons, we prove exact reduction results for two classes of graphon dynamical systems. We also show that, for finite-rank graphons, the reduced system preserves equilibria, linear stability, and local bifurcations. For infinite-rank graphons, we derive error bounds for the spectral reduction. Numerical examples using finite-rank and infinite-rank graphons illustrate the exact and approximate regimes of the theory.

\end{abstract}

\section{Introduction}\label{sec:intro}

The last three decades have seen extensive work on deterministic and stochastic dynamical systems on complex networks~\cite{boccaletti2006complex, barrat2008dynamical, Porter2016book}. A central motivation for such work is that both the dynamical rules and the underlying population structure, represented by networks and related objects such as hypergraphs and simplicial complexes, influence phenomena in many applications, including epidemic spreading, information cascades, and climate dynamics.

Dynamical systems on networks of practical interest are usually high-dimensional because, if each of the $N$ nodes carries at least one state variable, then the state space is already at least $N$-dimensional. Because real-world networks are often complex in multiple ways~\cite{Barabasi2016book, Newman2018book}, a recurring strategy for understanding dynamical systems on networks is dimension reduction, i.e., the construction of a much lower-dimensional dynamical system that approximates key qualitative and quantitative features of the full dynamics, such as equilibria, bifurcations, and resilience.
Such reductions are useful because they replace a large system of coupled equations by a small number of effective variables that can be analyzed, simulated, and interpreted more easily. For example, they can help locate bifurcations such as epidemic thresholds, estimate the resilience of equilibria, perform parameter sweeps, and compare the behavior of different networks without repeatedly solving the full high-dimensional system. Various dimension reduction methods of this kind have been examined for dynamics on networks~\cite{Golubitsky2006BullAmMathSoc, Gfeller2008PhysRevLett, Pecora2014NatCommun, Pastorsatorras2015RevModPhys, gao2016universal, laurence2019spectral, thibeault2020threefold, tu2017collapse, kundu2022accuracy, masuda2022dimension, vegue2023dimension, thibeault2024low}. Reduced-order modeling techniques for more general high-dimensional dynamical systems have also been studied extensively~\cite{rega2005dimension, Benner2015SIAMRev, Quarteroni2016book, rowley2017model, Brunton2022book, ghadami2022data, hesthaven2022reduced}.

A recent paper proposed a low-rank hypothesis for finite networks~\cite{thibeault2024low}, i.e., the matrix describing the interactions in a complex system, whether in models or empirical data, often has low effective rank. The authors further showed that this low-rank structure enables low-dimensional reduction of the dynamics. Their analysis rests on the singular value decomposition (SVD) of the interaction matrix of the network and, in particular, on the Schmidt--Eckart--Young--Mirsky theorem, which identifies the truncated SVD as the optimal low-rank approximation. We refer to this method as spectral reduction, following their earlier work~\cite{laurence2019spectral, thibeault2020threefold}; it is defined precisely in Section~\ref{sec:reduction}. Spectral reduction uses the dominant eigenvectors or singular vectors of the matrix coding the node-to-node interaction structure to construct reduced variables. Two results from~\cite{thibeault2024low} are especially relevant here. First, for a certain class of network dynamical systems, the dynamics can be exactly reduced to a system whose dimension is the rank of the interaction matrix. Second, for linear dynamics, the approximation error of the spectral reduction can be bounded explicitly in terms of the singular values of that matrix.
%

The purpose of this paper is to carry out this program rigorously for graphon dynamical systems. A graphon is a symmetric measurable function $W:\I^2\to\I$ that arises as the limit of a growing sequence of dense graphs~\cite{lovasz2006limits, borgs2008convergent, lovasz2012large}. We study spectral reduction for graphon dynamical systems, i.e., dynamical systems on graphons, for three reasons. First, because a graphon is symmetric, corresponding to an undirected network limit, the natural object replacing the adjacency matrix is a self-adjoint compact integral operator on $L^2(\I)$. In contrast, the low-rank hypothesis of~\cite{thibeault2024low} was formulated for general finite networks, including directed networks. We therefore expect the restriction to symmetric $W$ to simplify and sharpen the theory of spectral reduction. Second, passing from matrices to integral operators raises questions that do not arise in finite dimensions. Such questions include what the relevant eigenfunctions are, in what sense a graphon has an SVD or a rank, and how leading eigenvalues and eigenfunctions should be defined. These issues must be addressed carefully in the functional setting. Third, the graphon provides an $N$-independent description of the large-scale structure common to an entire convergent sequence of dense networks. Therefore, the reduction and its error bounds can be stated without reference to the number of nodes, $N$.
In this setting, we develop the spectral reduction directly at the graphon level and use the spectral structure of the graphon operator to identify when the reduction is exact and when its error is controlled by the decay of the graphon's spectrum.

\section{Graphon dynamical systems}\label{sec:setup}

We set up the graphon dynamical system to be reduced and state the standing assumptions under which the subsequent spectral reduction estimates are derived.
Let $I=\I$ and, for $p\ge1$, let $L^p(I)$ and $L^p(I^2)$ denote the usual Lebesgue spaces with norms $\norm{\cdot}_{L^p}$; we write $\inner{u}{v}=\int_I u(y)v(y)\,\dd y$ for the $L^2(I)$ inner product and $\norm{\cdot}=\norm{\cdot}_{L^2(I)}$ unless indicated otherwise. We denote by $\mathbf{1}$ the constant function $\mathbf{1}(y)=1$ on $I$.

In accordance with the standard convention for dense undirected graph limits, a graphon is a symmetric measurable function $W:I^2\to I$, i.e., $W(y,y')=W(y',y)$ and $0\le W(y,y')\le1$ for almost every $(y,y')\in I^2$.

The class of graphon dynamical systems that we primarily consider is
\begin{equation}\label{eq:graphon-ds}
\frac{\partial x(y,t)}{\partial t}=f\bigl(x(y,t)\bigr)+\int_0^1 W(y,y')G\bigl(x(y,t),x(y',t)\bigr)\,\dd y',
\qquad x(\cdot,0)=g\in S,
\end{equation}
where $x(y,t)\in\R$ is the state at location $y\in I$ and time $t$, $f:\Omega\to\R$ is the intrinsic, or self, dynamics, and $G:\Omega\times\Omega\to\R$ is the pairwise coupling, with $\Omega\subset\R$ a bounded interval. We impose the following assumptions~\cite{medvedev2014nonlinear, kaliuzhnyi2017semilinear, Eldo2025arxiv}.

\begin{assumption}[Lipschitz continuity]\label{as:lip}
The self-dynamics $f:\Omega\to\R$ is Lipschitz with constant $L_f$, i.e., $\abs{f(x)-f(x')}\le L_f\abs{x-x'}$ for all $x,x'\in\Omega$. The coupling $G:\Omega\times\Omega\to\R$ is Lipschitz with constant $L_G$ in the Euclidean norm, i.e.,
\begin{equation}
\abs{G(a,b)-G(a',b')}\le L_G\norm{(a,b)-(a',b')}_2
\end{equation}
for all $(a,b),(a',b')\in\Omega\times\Omega$. In particular, we will use the consequence
\begin{equation}
\abs{G(a,b)-G(a',b')}\le L_G(\abs{a-a'}+\abs{b-b'}).
\end{equation}
\end{assumption}

\begin{assumption}[Positively invariant domain]\label{as:inv}
There is a bounded interval $\Omega=[a_\ast,b_\ast]\subset\R$ such that the set
$S=\{u\in L^2(I): u(y)\in\Omega \text{ for a.e.\ }y\in I\}$ is positively invariant for~\eqref{eq:graphon-ds}. In other words, if $g\in S$, then the solution satisfies $x(\cdot,t)\in S$ for all $t\ge0$.
\end{assumption}

Under Assumptions~\ref{as:lip} and \ref{as:inv}, for every $g\in S$ and
every finite $T>0$, the initial value problem~\eqref{eq:graphon-ds}
has a unique solution
$x\in C^1\bigl([0,T];L^2(I)\bigr)$ with $x(\cdot,t)\in S$ for all $t\in[0,T]$, and the solution is uniformly bounded~\cite[Theorem~1 and Remark~1]{Eldo2025arxiv}, i.e.,
\begin{equation}
\sup_{t\in[0,T]}\norm{x(\cdot,t)}_{L^\infty(I)}\le\max\{\abs{a_\ast},\abs{b_\ast}\}=:B_\Omega.
\label{eq:B_Omega-def}
\end{equation}

\begin{remark}
Assumption~\ref{as:inv} is standing throughout, and is used only to keep solutions in the bounded set $S$ so that the Lipschitz constants $L_f$ and $L_G$ apply along trajectories. All the reduction estimates below are for a fixed solution $x(\cdot,t)\in S$ of~\eqref{eq:graphon-ds} on a fixed time horizon $[0,T]$.
\end{remark}

\section{Spectral theory of symmetric graphons}\label{sec:spectral-theory}

This section collects the operator-theoretic tools that the main developments, starting in Section~\ref{sec:reduction}, rely on. Each statement in this section is a specialization to graphons of a classical result on compact self-adjoint operators~\cite{conway2019course, reed1980methods, gohberg1978introduction, lovasz2012large}.
%
%
We give short proofs only to show the specialization to graphons and to keep the paper self-contained. We only include the facts used later in the paper.

\subsection{Graphon operator and its spectral decomposition}

A graphon enters the dynamics~\eqref{eq:graphon-ds} only through the integral of the coupling against $W$. It is therefore natural to package that operation as a single linear operator on the state space and to study the graphon through the properties of that operator, rather than through the function $W$ itself. Every reduction statement in this paper is phrased in terms of it.

\begin{definition}[Graphon operator~{\cite[Ch.~7]{lovasz2012large}}]\label{def:operator}
%
%
For a graphon $W$, the graphon operator $\Tw:L^2(I)\to L^2(I)$ is
\begin{equation}\label{eq:def-Tw}
(\Tw u)(y)=\int_0^1 W(y,y')u(y')\,\dd y',\qquad u\in L^2(I).
\end{equation}
\end{definition}

Three properties of $\Tw$, shown in the next lemma, do all of the work in what follows. Compactness gives a discrete nonzero spectrum that can accumulate only at the origin, so that there is a well-defined notion of leading eigenvalues at which to truncate. Self-adjointness, which owes to the symmetry of the graphon, makes the eigenvalues real and lets the eigenfunctions be chosen orthonormal, so that they can serve as a reduction basis. Being Hilbert--Schmidt makes the squared eigenvalues summable, which later converts a statement about the decay of the spectrum into a bound on the approximation error.

\begin{lemma}[{\src[\S VI.6]{reed1980methods}}]\label{lem:HS}
%
%
For any graphon $W$, the operator $\Tw$ is a self-adjoint Hilbert--Schmidt, hence compact, operator on $L^2(I)$, with Hilbert--Schmidt norm $\norm{\Tw}_{\mathrm{HS}}=\norm{W}_{L^2(I^2)}\le1$ and operator norm $\norm{\Tw}_{\mathrm{op}}\le\norm{W}_{L^2(I^2)}$.
\end{lemma}

\begin{proof}
Because $W$ is measurable and bounded by $1$ on the finite-measure set $I^2$, we have $W\in L^2(I^2)$ with $\norm{W}_{L^2(I^2)}\le1$. An integral operator with $L^2(I^2)$ kernel is Hilbert--Schmidt, with $\norm{\Tw}_{\mathrm{HS}}=\norm{W}_{L^2(I^2)}$~\cite[\S VI.6]{reed1980methods}.
%
%
In general, Hilbert--Schmidt operators are compact, and $\norm{\Tw}_{\mathrm{op}}\le\norm{\Tw}_{\mathrm{HS}}$. Self-adjointness follows from the symmetry $W(y,y')=W(y',y)$. Specifically, we obtain, for $u,v\in L^2(I)$,
\begin{align}
\inner{\Tw u}{v}
&=\iint W(y,y')u(y')v(y)\,\dd y'\,\dd y \notag\\
&=\iint W(y',y)u(y')v(y)\,\dd y\,\dd y' \notag\\
&=\inner{u}{\Tw v}. 
\end{align}
\end{proof}

Compactness and self-adjointness together yield the spectral theorem, the graphon counterpart of the diagonalization of a symmetric matrix by an orthonormal basis of eigenvectors. The spectral theorem produces the eigenfunctions $\varphi_k$ that serve as the reduction basis in Section~\ref{sec:reduction}, and the eigenvalues $\lambda_k$ whose rate of decay controls the reduction error in Section~\ref{sec:error}.

\begin{theorem}[Spectral decomposition of a graphon~{\cite[Ch.~II, \S5]{conway2019course}, \cite[Ch.~7]{lovasz2012large}}]\label{thm:spectral}
%
%
Let $W$ be a graphon. The nonzero eigenvalues of $\Tw$, counted with multiplicity, form either a finite sequence or a sequence tending to $0$. 
We denote them by $(\lambda_k)_{k\in K}$, where $K=\emptyset$ if
there are no nonzero eigenvalues, $K=\{1,\ldots,r\}$ if there are
exactly $r\in\{1,2,\ldots\}$ nonzero eigenvalues, counted with
multiplicity, and $K=\mathbb N$ if there are infinitely many eigenvalues.
They are ordered by nonincreasing modulus:
\begin{equation}
\abs{\lambda_1}\ge\abs{\lambda_2}\ge\cdots ,
\end{equation}
with the obvious finite interpretation when $K=\{1,\ldots,r\}$. There is an orthonormal system $(\varphi_k)_{k\in K}$ in $L^2(I)$ of corresponding eigenfunctions, $\Tw\varphi_k=\lambda_k\varphi_k$, such that
\begin{align}
W(y,y')&=\sum_{k\in K}\lambda_k\varphi_k(y)\varphi_k(y')
\quad\text{in } L^2(I^2), \label{eq:spectral-expansion}\\
\sum_{k\in K}\lambda_k^2&=\norm{W}_{L^2(I^2)}^2 . \label{eq:spectral-parseval}
\end{align}
Each nonzero eigenvalue has finite multiplicity, and
\begin{equation}
\overline{\operatorname{span}}\{\varphi_k:k\in K\}
=\overline{\ran\Tw}
=(\ker\Tw)^\perp .
\end{equation}
\end{theorem}

\begin{proof}
The operator $\Tw$ is compact and self-adjoint by
Lemma~\ref{lem:HS}. The spectral theorem for compact self-adjoint
operators yields real nonzero eigenvalues, each of finite
multiplicity and accumulating only at $0$, together with an
orthonormal eigenbasis of $(\ker\Tw)^\perp$
~\cite[Ch.~II, \S5]{conway2019course}. Ordering the nonzero
eigenvalues by decreasing modulus is a relabeling. For every
$u\in L^2(I)$,
\begin{equation}\label{eq:spectral-operator-expansion}
\Tw u
=
\sum_{k\in K}
\lambda_k\inner{\varphi_k}{u}\varphi_k,
\end{equation}
with convergence in $L^2(I)$. Because $\Tw$ is Hilbert--Schmidt,
\begin{equation}\label{eq:HS-spectral-sum}
\norm{\Tw}_{\mathrm{HS}}^2
=
\sum_{k\in K}\lambda_k^2.
\end{equation}
The Hilbert--Schmidt norm of an integral operator equals the
$L^2(I^2)$ norm of its kernel
~\cite[\S VI.6]{reed1980methods}.
%
%
Consequently, \eqref{eq:spectral-operator-expansion} gives the kernel
expansion~\eqref{eq:spectral-expansion}, and
\eqref{eq:HS-spectral-sum} gives
\eqref{eq:spectral-parseval}. The characterization of the closed span
of the nonzero eigenfunctions follows from the spectral theorem for compact self-adjoint operators.
\end{proof}

We fix an orthonormal basis in each nonzero eigenspace of $\Tw$ and a
rule for ordering eigenfunctions whose eigenvalues have the same
modulus. The ordered sequence $(\varphi_k)$, and every subsequent
truncation based on its first $n$ elements, are understood relative
to these choices. If a cutoff splits a collection of eigenfunctions
whose eigenvalues have the same modulus, the retained subspace and the
resulting reduced system may depend on the chosen basis or ordering.
When the entire tied collection of eigenfunctions is retained, changing its orthonormal
basis only changes the coordinates used for representing the same retained
subspace.

If $K=\emptyset$, equivalently if $W=0$ a.e., we set
$\lambda_k=0$ for all $k\ge1$ and choose a fixed orthonormal basis
$(\varphi_k)_{k\ge1}$ of $L^2(I)$. If
$K=\{1,\ldots,r\}$, we use the convention $\lambda_k=0$ for $k>r$
whenever formulas such as $\lambda_{n+1}$ are used, and we append a
fixed orthonormal basis of $\ker\Tw$ after the $r$ nonzero
eigenfunctions. All statements and estimates below remain valid for
any such fixed choice, although the projection onto the first $n$
eigenfunctions and quantities defined from it may depend on the
choice. When $K=\{1,\ldots,r\}$, the subspace spanned by the first
$r$ eigenfunctions is $\ran\Tw$, whereas any additional zero
eigenfunctions used when $n>r$ are basis-dependent.

We next record the approximation error obtained by truncating the
spectral expansion.

\begin{corollary}[Spectral truncation error]
\label{cor:spectral-truncation}
For $n\ge1$, define
\begin{equation}\label{eq:spectral-truncation-kernel}
W_n
:=
\sum_{\substack{k\in K\\ k\le n}}
\lambda_k\varphi_k\otimes\varphi_k,
\end{equation}
where
$(\varphi\otimes\psi)(y,y')
:=
\varphi(y)\psi(y')$.
Then,
\begin{equation}\label{eq:L2-convergence}
\norm{W-W_n}_{L^2(I^2)}^2
=
\int_0^1\!\!\int_0^1
\bigl\lvert W(y,y')-W_n(y,y')\bigr\rvert^2
\,\dd y\,\dd y'
=
\sum_{\substack{k\in K\\ k>n}}\lambda_k^2
\longrightarrow0
\qquad
(n\to\infty),
\end{equation}
where an empty tail is interpreted as $0$.
\end{corollary}

\begin{proof}
The kernels
$\varphi_k\otimes\varphi_k$, $k\in K$, are orthonormal in
$L^2(I^2)$. Therefore, the second equality in
\eqref{eq:L2-convergence} follows from the kernel expansion
\eqref{eq:spectral-expansion} and the Pythagorean theorem. The
convergence to zero follows from
\eqref{eq:spectral-parseval}.
\end{proof}

We call $\lambda_1$ an eigenvalue of largest modulus and $\varphi_1$
a leading eigenfunction when $W\ne0$. If several eigenfunctions have the same
largest eigenvalue modulus, the leading eigenfunction depends on the
basis and ordering fixed above.

We also record the following consequence of the spectral decomposition, which is used in the reduction error estimates.

\begin{corollary}[Restricted operator norm]\label{cor:restricted-norm}
Let $\Pn$ be the orthogonal projection onto $V_n:=\operatorname{span}\{\varphi_1,\ldots,\varphi_n\}$, with the finite-rank convention of Theorem~\ref{thm:spectral}. Then,
\begin{equation}\label{eq:restricted-norm}
\norm{\Tw(I-\Pn)}_{\mathrm{op}}
=
\norm{(I-\Pn)\Tw}_{\mathrm{op}}
=
\abs{\lambda_{n+1}},
\end{equation}
where $\lambda_{n+1}=0$ when no nonzero eigenfunction remains after the first
$n$ eigenfunctions.
\end{corollary}

\begin{proof}
For $u\in L^2(I)$,
\begin{equation}
\Tw(I-\Pn)u=\sum_{k>n}\lambda_k\inner{\varphi_k}{u}\varphi_k,
\label{eq:restricted-norm-proof-1}
\end{equation}
with an empty sum when no nonzero eigenfunction remains after the first $n$
eigenfunctions.
Equation~\eqref{eq:restricted-norm-proof-1} leads to
\begin{equation}
\norm{\Tw(I-\Pn)u}^2
=
\sum_{k>n}\lambda_k^2\inner{\varphi_k}{u}^2
\le
\lambda_{n+1}^2\norm{u}^2.
\end{equation}
Equality holds by taking $u=\varphi_{n+1}$ when a nonzero $(n+1)$st mode exists. If no such mode exists, the operator is zero. Moreover, $\Tw$ and $\Pn$ commute because both are diagonal in the eigenfunction basis. Hence $(I-\Pn)\Tw=\Tw(I-\Pn)$, and the same norm identity holds for $(I-\Pn)\Tw$.
\end{proof}

\subsection{Singular values and the SVD of a graphon}

For a general non-symmetric kernel $W\in L^2(I^2)$, $\Tw$ is compact
but not self-adjoint, and one works with its SVD: there are
orthonormal systems $(u_k)$ and $(v_k)$ and singular values
$\sigma_k\ge0$, $\sigma_k\downarrow0$, such that, for every
$z\in L^2(I)$,
$\Tw z
=
\sum_k\sigma_k\inner{v_k}{z}u_k$;
the $\sigma_k$ are the eigenvalues of $(\Tw^\ast\Tw)^{1/2}$.
This is the functional analogue of the matrix SVD used in~\cite{thibeault2024low}. For a graphon, which is symmetric by definition, the SVD carries no information beyond the eigendecomposition, except for the signs of the nonzero eigenvalues.

\begin{definition}[Singular values of a graphon~{\cite[Ch.~II]{gohberg1978introduction}}]\label{def:sv}
%
%
The singular values of a graphon $W$ are the eigenvalues of $(\Tw^\ast\Tw)^{1/2}$, arranged in nonincreasing order $\sigma_1\ge\sigma_2\ge\cdots\ge0$ and counted with multiplicity.
\end{definition}

We can now state the precise relation between the SVD and the eigendecomposition for graphons.

\begin{proposition}[SVD equals eigendecomposition for symmetric $W$]\label{prop:svd}
Let $W$ be a graphon with nonzero eigenvalue--eigenfunction pairs $(\lambda_k,\varphi_k)_{k\in K}$ as in Theorem~\ref{thm:spectral}. Then, the nonzero singular values of $W$ are
\begin{equation}
\sigma_k=\abs{\lambda_k},\qquad k\in K,
\end{equation}
ordered with multiplicity. For every nonzero eigenvalue, $\Tw$ has the singular representation
\begin{equation}\label{eq:svd}
\Tw u=\sum_{k\in K}\abs{\lambda_k}\inner{v_k}{u}\varphi_k,
\qquad
v_k:=\operatorname{sign}(\lambda_k)\varphi_k .
\end{equation}
Therefore, the left and right singular functions associated with nonzero singular values differ only by signs, and they coincide whenever $\lambda_k>0$. The zero singular subspace is $\ker\Tw$ and may be represented by any orthonormal basis of that kernel.
\end{proposition}

\begin{proof}
Since $\Tw$ is self-adjoint, $\Tw^\ast\Tw=\Tw^2$. Therefore, if $\Tw\varphi_k=\lambda_k\varphi_k$, then $\Tw^\ast\Tw\varphi_k=\lambda_k^2\varphi_k$, and $(\Tw^\ast\Tw)^{1/2}$ has eigenvalue $\abs{\lambda_k}$ on the same eigenspace. So, the nonzero singular values are $\abs{\lambda_k}$. 
For $\lambda_k\ne0$, we write
$\lambda_k=\operatorname{sign}(\lambda_k)\abs{\lambda_k}$ and set
$v_k=\operatorname{sign}(\lambda_k)\varphi_k$. Then,
\begin{align}
\Tw u
&=\sum_{k\in K}\lambda_k\inner{\varphi_k}{u}\varphi_k \notag\\
&=\sum_{k\in K}\abs{\lambda_k}\inner{v_k}{u}\varphi_k.
\end{align}
This is the singular-value expansion with left singular functions
$\varphi_k$, right singular functions $v_k$, and singular values
$\abs{\lambda_k}$. The left singular functions $(\varphi_k)$ are orthonormal by
the spectral theorem. The right singular functions $(v_k)$ are also
orthonormal because multiplying each $\varphi_k$ by 
$\operatorname{sign}(\lambda_k) \in \{ -1, 1 \}$ preserves orthonormality.
\end{proof}

Proposition~\ref{prop:svd} explains why, for graphons, we may work throughout with the real eigenpairs $(\lambda_k,\varphi_k)$ and use $\sigma_k=\abs{\lambda_k}$ for the corresponding singular values. For a general matrix, the left and right singular vectors differ, which is one source of approximation error in~\cite{thibeault2024low}.

\subsection{Rank of a graphon and finite-rank graphons}

The low-rank hypothesis is a statement about the rank of finite networks. Therefore, we need to say what the rank of a graphon is. For a matrix, the dimension of the range, the number of nonzero eigenvalues, and the number of nonzero singular values all agree. The same equivalence holds for compact self-adjoint graphon operators.

\begin{definition}[Rank]\label{def:rank}
The rank of a graphon $W$ is
\begin{equation}
\rank(W)=\dim\ran\Tw\in\{0,1,2,\ldots\}\cup\{\infty\},
\end{equation}
equivalently, the number of nonzero eigenvalues of $\Tw$ counted with multiplicity, or the number of nonzero singular values.
\end{definition}

Because the rank of a graphon can be infinite, the low-rank hypothesis corresponds in the graphon setting either to finite-rank graphons or to infinite-rank graphons with rapidly decaying singular values. In this article, finite-rank graphons provide the setting in which exact reduction is possible.

\begin{proposition}[Finite-rank graphons]\label{prop:finite-rank}
For a graphon $W$ and $r\in\{1,2,\ldots\}$, the following are equivalent.
\begin{enumerate}
\item[\textup{(i)}] $\rank(W)=r$.
\item[\textup{(ii)}] $W(y,y')=\sum_{k=1}^{r}\lambda_k\varphi_k(y)\varphi_k(y')$ a.e., with $\lambda_k\neq0$ and $(\varphi_k)_{k=1}^r$ orthonormal.
\end{enumerate}
In this case, $\ran\Tw=\operatorname{span}\{\varphi_1,\ldots,\varphi_r\}$ is $r$-dimensional and closed. Moreover, if $n\ge r$ and $\Pn$ is the orthogonal projection onto a subspace containing $\operatorname{span}\{\varphi_1,\ldots,\varphi_r\}$, then $x-\Pn x\in\ker\Tw$ for every $x\in L^2(I)$.
\end{proposition}

\begin{proof}
The equivalence of \textup{(i)} and \textup{(ii)} follows directly from Theorem~\ref{thm:spectral}. If $\rank(W)=r$, then $\Tw$ has exactly $r$ nonzero eigenvalues, counted with multiplicity, and the spectral expansion terminates after $r$ terms. Conversely, if $W$ has the representation in \textup{(ii)}, then the range of $\Tw$ is contained in $\operatorname{span}\{\varphi_1,\ldots,\varphi_r\}$, and 
$\lambda_k\neq0$ for $k \in \{ 1,\ldots,r \}$. Therefore, $\rank(W)=r$.

The final statement holds because, for finite-rank $W$, $\ran\Tw$ is finite-dimensional, hence closed, and Theorem~\ref{thm:spectral} gives
$\ran\Tw=(\ker\Tw)^\perp=\operatorname{span}\{\varphi_1,\ldots,\varphi_r\}$.
If $n\ge r$ and $\Pn$ projects onto any subspace containing $\ran\Tw$, then $(I-\Pn)x\perp\ran\Tw$, so $(I-\Pn)x\in(\ran\Tw)^\perp=\ker\Tw$.
\end{proof}

Proposition~\ref{prop:finite-rank} is the graphon counterpart of the elementary fact that a symmetric matrix has finite rank precisely when its spectral expansion has only finitely many nonzero eigenvalues.

\section{Spectral dimension reduction}\label{sec:reduction}

This section sets up the reduction itself. We choose a small number of
leading eigenfunctions of the graphon and use them as weight
functions, so that each of them extracts one linear observable of the
state. These observables are the variables of the reduced
description. We derive the closed system of ordinary differential
equations that we consider as the dimension-reduced dynamical system.

The formulas below evaluate $f$ and $G$ at projected or reconstructed states. These reconstructed states need not take values in the interval $\Omega$, even when the true graphon state does. 
To make the reduced vector fields well defined on all of $\R^n$, we
choose Lipschitz extensions of $f$ to $\R$ and of $G$ to $\R^2$ that
preserve the Lipschitz constants $L_f$ and $L_G$, respectively, and
keep the same notation for the extensions.
These extensions do not change the original graphon dynamics on the positively invariant set $S$.

\subsection{Reduction operator and reduced system}

Fix a reduced dimension $n\ge1$ and the leading eigenfunctions $\varphi_1,\ldots,\varphi_n$ of $W$. Define the reduction operator $\Mn:L^2(I)\to\R^n$ and its adjoint $\Mn^\ast:\R^n\to L^2(I)$ by
\begin{equation}\label{eq:reduction-op}
(\Mn u)_\mu=\inner{\varphi_\mu}{u}=\int_0^1\varphi_\mu(y)u(y)\,\dd y,\qquad
\mu \in \{1,\ldots,n \},
\end{equation}
and
\begin{equation}
\Mn^\ast \xi=\sum_{\mu=1}^{n}\xi_\mu\varphi_\mu,
\label{eq:reduction-op-2}
\end{equation}
where $u\in L^2(I)$ and $\xi=(\xi_1,\ldots,\xi_n)^\top\in\R^n$ are arbitrary.
Because $(\varphi_\mu)$ is orthonormal, $\Mn\Mn^\ast=I_n$ (the $n\times n$ identity matrix) and $\Mn^\ast\Mn=\Pn$, the orthogonal projection onto $V_n=\operatorname{span}\{\varphi_1,\ldots,\varphi_n\}$. Therefore, $\Mn^\ast$ is the Moore--Penrose pseudoinverse of $\Mn$, the functional analogue of the reduction matrix and its pseudoinverse in~\cite{thibeault2024low}.

\begin{definition}[Observable]\label{def:observable}
An observable of the graphon dynamical system~\eqref{eq:graphon-ds} is a scalar quantity, computed from the state $x(\cdot,t)$, that is used to summarize that state by a single number evolving in time. All observables in this paper are linear, that is, of the form $x\mapsto\inner{\psi}{x}=\int_0^1\psi(y)x(y)\,\dd y$ for a fixed weight function $\psi\in L^2(I)$. Therefore, each observable is a weighted integral of the state.
The $n$ observables of the spectral reduction are those obtained by taking $\psi=\varphi_\mu$, the $n$ leading eigenfunctions of $W$:
\begin{equation}\label{eq:observables}
X_\mu(t)=\inner{\varphi_\mu}{x(\cdot,t)}=\int_0^1\varphi_\mu(y)x(y,t)\,\dd y,\qquad \mu=1,\ldots,n .
\end{equation}
\end{definition}

We denote by $X(t)=(X_1(t),\ldots,X_n(t))^\top=\Mn x(\cdot,t)\in\R^n$ the vector collecting the $n$ observables, and by $\widehat x:=\Mn^\ast X=\Pn x \in L^2(I)$ the reconstruction of the state from them.

Using the Lipschitz extensions fixed above, we define the extended
vector field $\mathcal N:L^2(I)\to L^2(I)$ by
\begin{equation}\label{eq:field}
\mathcal N(x)(y)
=
f\bigl(x(y)\bigr)
+
\int_0^1
W(y,y')G\bigl(x(y),x(y')\bigr)\,\dd y'.
\end{equation}
To see that this map is well defined, we let
$s(y):=(\Tw\mathbf{1})(y)$ and
$s_\ast:=\operatorname{ess\,sup}_{y\in I}s(y)\le1$.
Then, we obtain
\begin{equation}
\mathcal N(0)=f(0)\mathbf{1}+G(0,0)s\in L^2(I),
\end{equation}
and, for all $u,v\in L^2(I)$,
\begin{equation}\label{eq:field-lipschitz}
\norm{\mathcal N(u)-\mathcal N(v)}_{L^2}
\le
\left[
L_f+L_G\bigl(s_\ast+\norm{\Tw}_{\mathrm{op}}\bigr)
\right]
\norm{u-v}_{L^2}.
\end{equation}
Thus, $\mathcal N$ is globally Lipschitz on $L^2(I)$. On $S$, this
extended vector field agrees with the original right-hand side
of~\eqref{eq:graphon-ds}; hence, the original graphon dynamics is the
restriction of
$\partial_t x=\mathcal N(x)$
to the positively invariant set $S$.
The exact evolution of the observables is
\begin{equation}\label{eq:exact-obs}
\dot X_\mu=\inner{\varphi_\mu}{\mathcal N(x)},\qquad\mu=1,\ldots,n,
\end{equation}
which is not closed in $X$ because $\mathcal N(x)$ depends on the full state $x$. The $n$-dimensional spectral reduction closes~\eqref{eq:exact-obs} by replacing $x$ with its reconstruction $\Mn^\ast X$:
\begin{align}
\dot X_\mu &=\bigl(\Mn\mathcal N(\Mn^\ast X)\bigr)_\mu \notag\\
&= \int_0^1\varphi_\mu(y)\left[f\!\Bigl(\textstyle\sum_{\nu=1}^{n}X_\nu\varphi_\nu(y)\Bigr)
+\int_0^1 W(y,y')G\!\Bigl(\textstyle\sum_{\nu=1}^{n}X_\nu\varphi_\nu(y),\sum_{\nu=1}^{n}X_\nu\varphi_\nu(y')\Bigr)\dd y'\right]\dd y .
\label{eq:reduced}
\end{align}
System~\eqref{eq:reduced} is a closed system of $n$ ordinary differential equations for $X=(X_1,\ldots,X_n)$, which we refer to as the spectral reduction, following the terminology of~\cite{laurence2019spectral, thibeault2020threefold} for finite networks. Equivalently, in operator form
\begin{equation}\label{eq:reduced-op}
\dot X=F_n(X):=\Mn\mathcal N(\Mn^\ast X).
\end{equation}
This is the least-squares optimal closure in the following sense, which is the graphon analogue of~\cite[Suppl.\ Thm.\ 52]{thibeault2024low}.

\begin{proposition}[Optimal closure]\label{prop:closure}
Among all vector fields $F:\R^n\to\R^n$, the choice $F=F_n$ in~\eqref{eq:reduced-op} minimizes
\begin{equation}
\norm{\mathcal N(\Pn x)-\Mn^\ast F(\Mn x)}_{L^2(I)}
\end{equation}
for every $x$. More precisely, this choice makes the $V_n$-component of the closure error vanish:
\begin{equation}
\Mn^\ast F_n(\Mn x)=\Pn\mathcal N(\Pn x).
\end{equation}
\end{proposition}

\begin{proof}
For $u\in L^2(I)$, $\Mn^\ast F(\Mn u)\in V_n$, so
\begin{equation}
\norm{\mathcal N(\Pn u)-\Mn^\ast F(\Mn u)}^2=\norm{(I-\Pn)\mathcal N(\Pn u)}^2+\norm{\Pn\mathcal N(\Pn u)-\Mn^\ast F(\Mn u)}^2
\label{eq:optimal-closure-1}
\end{equation}
by orthogonality. The first term on the right-hand side of \eqref{eq:optimal-closure-1}
is independent of $F$. The second term is minimized to $0$ by $F=F_n$, which one can verify by substituting
$F_n(X)= \Mn\mathcal N(\Mn^\ast X)$ from \eqref{eq:reduced-op} and $\Pn = \Mn^\ast\Mn$ in the second term on the right-hand side of \eqref{eq:optimal-closure-1}.

\end{proof}

\section{Exact reduction of low-rank graphon dynamics}\label{sec:exact}

\subsection{Exact finite-dimensional reduction}\label{sec:exact-finite-dimensional}

We now prove that two natural classes of graphon dynamical systems reduce exactly to finite-dimensional dynamical systems when $W$ is of finite rank. This is the graphon analogue of the exact rank reduction of~\cite{thibeault2024low} and vindicates the low-rank hypothesis.

\begin{definition}[Reducible vector fields]\label{def:reducible}
We consider two classes of graphon dynamical systems.
\begin{itemize}
\item[\textup{(A)}] Neighbor-coupling form with linear self-dynamics:
\begin{equation}
\partial_t x(y,t)=cx(y,t)+D\int_0^1 W(y,y')h\left(x(y',t)\right)\,\dd y'
= cx+D\Tw\left(h\circ x\right),
\label{eq:reducible-A}
\end{equation}
where $c,D\in\R$ and $h:\R\to\R$ is Lipschitz.
\item[\textup{(B)}] Input-coupling form:
\begin{equation}
\partial_t x(y,t)=f\left(x(y,t)\right)+\tilde{g}\left(y,(\Tw x)(y,t)\right),
\label{eq:reducible-B}
\end{equation}
where $f$ is affine, $f(x)=cx+c_0$, and $\tilde{g}:I\times\R\to\R$ is measurable, satisfies $\tilde g(\cdot,0)\in L^2(I)$, and is Lipschitz in its second argument with constant $L_{\tilde{g}}$.
\end{itemize}
\end{definition}

Class \textup{(A)} is the subclass of~\eqref{eq:graphon-ds} in which $G(x(y),x(y'))=h(x(y'))$ depends only on the neighbor's state and the self-dynamics is linear.

Class \textup{(B)} is the graphon counterpart of the exactly
reducible finite-network class considered in~Eq.~(22) of \cite{thibeault2024low}. Class \textup{(B)} is in general not a subclass of~\eqref{eq:graphon-ds}. When $\tilde{g}$ is nonlinear in its second argument, the nonlinearity acts on the aggregated input $(\Tw x)(y)$ after the integration over $y'$ has been carried out. Such a term cannot be written as an integral $\int_0^1W(y,y')G(x(y),x(y'))\dd y'$ of a pointwise coupling.
A simple special case common to both classes is $c_0=0$ and
$\tilde g(y,z)=Dz$. This is class \textup{(A)} with $h(x)=x$ and is
also of the form~\eqref{eq:graphon-ds} with $G(a,b)=Db$.

Both classes are globally well posed on $L^2(I)$. We write
$\mathcal N_A$ and $\mathcal N_B$ for their vector fields. We first
note that $\mathcal N_A(0)=D h(0)\Tw\mathbf{1}\in L^2(I)$ and $\mathcal N_B(0)=c_0\mathbf{1}+\tilde g(\cdot,0)\in L^2(I)$.
Moreover, if $L_h$ is a Lipschitz constant for $h$, then, for all
$u,v\in L^2(I)$, we obtain
\begin{equation}
\norm{\mathcal N_A(u)-\mathcal N_A(v)}_{L^2}
\le
\left(
\abs{c}
+
\abs{D}L_h\norm{\Tw}_{\mathrm{op}}
\right)
\norm{u-v}_{L^2},
\end{equation}
and
\begin{equation}
\norm{\mathcal N_B(u)-\mathcal N_B(v)}_{L^2}
\le
\left(
\abs{c}
+
L_{\tilde g}\norm{\Tw}_{\mathrm{op}}
\right)
\norm{u-v}_{L^2}.
\end{equation}
Taking $v=0$ shows that both vector fields map $L^2(I)$ into itself,
and the same estimates show that they are globally Lipschitz. The
Picard--Lindel\"of theorem for ordinary differential equations in
Banach spaces therefore gives a unique global solution
$x\in C^1\bigl([0,\infty);L^2(I)\bigr)$
for every initial condition $x(\cdot,0)\in L^2(I)$.

We now state the exact reduction theorem. The mechanism differs between the two classes. In class \textup{(A)}, the coupling term always lands in $\ran\Tw=V_r$, so $V_r$ is an invariant subspace of the dynamics. The dynamics in $V_r$ is genuinely $r$-dimensional in this case. In class \textup{(B)}, the state leaves $V_r$, and no invariant subspace is available. In this case, what rescues the reduction is that the coupling sees the state only through $\Tw x$, which is determined by the $r$ observables. The observables then close among themselves even though the state does not.

\begin{theorem}[Exact finite-dimensional reduction]\label{thm:exact}
Let $W$ have finite rank $r$, with nonzero eigenpairs $(\lambda_k,\varphi_k)_{k=1}^{r}$, and let $V_r=\operatorname{span}\{\varphi_1,\ldots,\varphi_r\}$. Consider the observables $X_\mu(t)=\inner{\varphi_\mu}{x(\cdot,t)}$, $\mu=1,\ldots,r$.
\begin{enumerate}
\item[\textup{(i)}] For class \textup{(A)}, if the initial condition satisfies $x(\cdot,0) \in V_r$, then $x(\cdot,t)\in V_r$ for all $t\ge0$ (the subspace $V_r$ is positively invariant), the state is exactly reconstructed by $x=\sum_{\nu=1}^r X_\nu\varphi_\nu$, and
\begin{equation}\label{eq:exactA}
\dot X_\mu=cX_\mu+D\lambda_\mu\int_0^1\varphi_\mu(y)h\!\Bigl(\textstyle\sum_{\nu=1}^{r}X_\nu\varphi_\nu(y)\Bigr)\dd y,\qquad\mu=1,\ldots,r.
\end{equation}
\item[\textup{(ii)}] For class \textup{(B)}, the observables satisfy the closed, autonomous $r$-dimensional system
\begin{equation}\label{eq:exactB}
\dot X_\mu=cX_\mu+c_0\inner{\varphi_\mu}{\mathbf{1}}+\int_0^1\varphi_\mu(y)
\tilde{g}\!\Bigl(y,\ \textstyle\sum_{\nu=1}^{r}\lambda_\nu X_\nu\varphi_\nu(y)\Bigr)\dd y,\qquad\mu=1,\ldots,r,
\end{equation}
for every initial condition $x(\cdot,0)\in L^2(I)$.
\end{enumerate}
\end{theorem}

\begin{proof}
Because $\rank(W)=r$, Proposition~\ref{prop:finite-rank} gives $\ran\Tw=V_r$ and
\begin{equation}
\Tw x=\sum_{\nu=1}^{r}\lambda_\nu\inner{\varphi_\nu}{x}\varphi_\nu=\sum_{\nu=1}^{r}\lambda_\nu X_\nu\varphi_\nu
\label{eq:Twx-2}
\end{equation}
for every $x\in L^2(I)$. We note that $\Tw x$ depends on $x$ only through $X_1,\ldots,X_r$.

(i) If $x\in V_r$, then $cx\in V_r$, and $\Tw(h\circ x)\in\ran\Tw=V_r$. 
Therefore, the vector field maps $V_r$ into itself. By the global
well-posedness established above and uniqueness of solutions, a
solution with $x(\cdot,0)\in V_r$ remains in $V_r$ for all $t\ge0$.
Thus, we obtain
$x=\sum_{\nu=1}^r X_\nu\varphi_\nu$.
By differentiating the observable $X_\mu(t)=\inner{\varphi_\mu}{x(\cdot,t)}$ with respect to $t$ and using $\partial_t x=cx+D \Tw(h\circ x)$, we obtain
\begin{align}
\dot X_\mu &=\inner{\varphi_\mu}{cx+D\Tw(h\circ x)} \notag\\
&= c X_\mu+D\inner{\Tw\varphi_\mu}{h\circ x} \notag\\
&= c X_\mu+D\lambda_\mu\inner{\varphi_\mu}{h\circ x},
\label{eq:(i)-proof-1}
\end{align}
using self-adjointness $\inner{\varphi_\mu}{\Tw(h\circ x)}=\inner{\Tw\varphi_\mu}{h\circ x}$ and $\Tw\varphi_\mu=\lambda_\mu\varphi_\mu$. Substituting $x=\sum_{\nu=1}^{r}X_\nu\varphi_\nu$ into the right-hand side of~\eqref{eq:(i)-proof-1} yields~\eqref{eq:exactA}.

(ii) By differentiating the observable and using
$\partial_t x=cx+c_0+\tilde{g}(\cdot,\Tw x)$, we obtain
\begin{align}
\dot X_\mu &= \inner{\varphi_\mu}{\partial_t x}\notag\\
&= c\inner{\varphi_\mu}{x}+c_0\inner{\varphi_\mu}{\mathbf{1}}+\inner{\varphi_\mu}{\tilde{g}(\cdot,\Tw x)}\notag\\
&= cX_\mu+c_0\inner{\varphi_\mu}{\mathbf{1}}+\int_0^1\varphi_\mu(y)\tilde{g}\bigl(y,(\Tw x)(y)\bigr)\dd y .
\label{eq:(ii)-proof-1}
\end{align}
Substituting $(\Tw x)(y)=\sum_{\nu = 1}^r \lambda_\nu X_\nu\varphi_\nu(y)$ in the right-hand side of \eqref{eq:(ii)-proof-1} gives~\eqref{eq:exactB}. The right-hand side depends only on $X$, so the system is closed and autonomous.

\end{proof}

\begin{remark}[Reconstruction and the kernel of $\Tw$]\label{rem:kernel}
In case \textup{(ii)}, the observables $X$ evolve autonomously and exactly, but they do not in general determine the full state. The component $x_\perp=(I-P_r)x\in\ker\Tw$ obeys
\begin{equation}
\partial_t x_\perp=cx_\perp+c_0(I-P_r) \mathbf{1}+(I-P_r)\tilde{g}(\cdot,\Tw x).
\label{eq:remark-rem-kernel}
\end{equation}
Using \eqref{eq:Twx-2}, we write the solution of
\eqref{eq:remark-rem-kernel} as
\begin{align}
x_\perp(t)
&=e^{ct}x_\perp(0)
+\int_0^t e^{c(t-\tau)}
\biggl[
c_0(I-P_r)\mathbf{1}
+(I-P_r)\tilde{g}\left(\cdot,\sum_{\nu=1}^{r}\lambda_\nu X_\nu(\tau)\varphi_\nu\right)
\biggr]\,\dd\tau ,
\label{eq:xperp-variation}
\end{align}
which shows precisely what information is missing from the observables. The integral term is determined by the trajectory $X(\tau)$, but the homogeneous term $e^{ct}x_\perp(0)$ depends on the initial kernel component $x_\perp(0)$, which is invisible to the observables. For example, if $u \in\ker\Tw$, then $x(0)$ and $x(0)+u$ have the same values of $X_\mu(0)=\inner{\varphi_\mu}{x(0)}$, $\mu=1,\ldots,r$, and hence generate the same observable trajectory $X(t)$. However, their full states differ by $e^{ct}u$ in the kernel direction.

Thus, in class \textup{(B)}, the trajectory of the observables, but not the observables at an instantaneous time $t$, determines the full state at time $t$ only after the initial kernel component has also been specified. If one imposes the additional condition $x_\perp(0)=0$, then $x_\perp(t)$ is determined by $X(\tau)$ through~\eqref{eq:xperp-variation}, and the full state can be reconstructed from the observable trajectory. If $x_\perp(0)\neq0$ and $c<0$, the contribution
$e^{ct}x_\perp(0)$ decays to
zero as $t\to\infty$. Therefore, the full state is asymptotically
reconstructible from the observable trajectory.
\end{remark}

The most important special case is linear dynamics, for which the
equations for the retained observables are exact for every
truncation dimension $n$.

\begin{corollary}[Exact spectral decoupling of linear graphon dynamics]\label{cor:linear-exact}
Consider the linear graphon dynamical system
\begin{equation}\label{eq:linear}
\partial_t x(y,t)=cx(y,t)+D\int_0^1 W(y,y')x(y',t)\,\dd y'=cx+D\Tw x
\end{equation}
for any graphon $W$ (finite or infinite rank). Then, every observable $X_{\mu}(t)=\inner{\varphi_{\mu}}{x(\cdot,t)}$, where $\mu \ge 1$, evolves independently and exactly according to
\begin{equation}\label{eq:linear-modes}
\dot X_{\mu}=(c+D\lambda_{\mu})X_{\mu},
\end{equation}
namely,
\begin{equation}
X_{\mu}(t)=X_{\mu}(0)e^{(c+D\lambda_{\mu})t}.
\end{equation}
In addition, the component of $x$ in $\ker\Tw$ evolves as $\partial_t x_\perp=cx_\perp$. 
As a result, for every $n$, the $n$-dimensional spectral
reduction~\eqref{eq:reduced} reproduces the dynamics of the first $n$
observables exactly.
\end{corollary}

\begin{proof}
$\dot X_{\mu} = \inner{\varphi_{\mu}}{cx+D\Tw x}=cX_{\mu}+D\inner{\Tw\varphi_{\mu}}{x}=(c+D\lambda_{\mu})X_{\mu}$, giving~\eqref{eq:linear-modes}.
\end{proof}

Corollary~\ref{cor:linear-exact} is the self-adjoint graphon
counterpart of the finite-network linear result
of~\cite{thibeault2024low}. The symmetry assumed here makes the result particularly transparent, whereas the finite-network result~\cite{thibeault2024low}
also permits asymmetric interaction matrices, i.e., directed networks.

\subsection{Equilibria, stability, and bifurcations}\label{sec:bifurcation}

Theorem~\ref{thm:exact} guarantees that the reduced system exactly reproduces the observables. However, by itself, this result does not imply that the reduced system is a faithful guide to the qualitative behavior of the full graphon dynamics. It is the qualitative behavior (e.g., which equilibria exist, whether they are stable, and at which parameter values they lose stability) that a dimension reduction often aims to predict~\cite{gao2016universal, laurence2019spectral, liu2022network}. The gap is genuine. For example, an equilibrium could be stable as seen by the $r$ observables and yet unstable in the state space $L^2(I)$.

We show that, for the exactly reducible classes, such a mismatch between the stability of the reduced dynamical system and that of the full graphon dynamical system does not occur in any nontrivial way. The reason is structural: the linearization of the graphon vector field at an equilibrium is block triangular with respect to the splitting $L^2(I)=V_r\oplus V_r^\perp$, and its $V_r^\perp$ block is exactly $cI$. The discarded directions therefore have the single exponential rate
$c$, which is independent of the equilibrium. We isolate this triangular structure first.

\begin{lemma}[Triangular structure of the linearization]\label{lem:triangular}
Let $W$ have finite rank $r$, so that $L^2(I)=V_r\oplus V_r^\perp$ with $V_r=\ran\Tw$ and $V_r^\perp=\ker\Tw$.
%
%
Regard each vector field in Definition~\ref{def:reducible} (i.e., \eqref{eq:reducible-A} and \eqref{eq:reducible-B}) as a map
from $L^2(I)$ to $L^2(I)$. Let
$x^\ast\in S$ be an equilibrium of the graphon dynamics.
\begin{enumerate}
	\item[\textup{(i)}]
	For class \textup{(A)}, suppose that $h\in C^1$ on a neighborhood
	of $\Omega$ and that $x^\ast\in V_r$. Then, the vector field
$\mathcal N_A(x)=cx+D\Tw(h\circ x)$
	is Fr\'echet differentiable at $x^\ast$, and its derivative is the
	bounded operator
	\begin{equation}
		\mathcal L v
		:=\mathrm{D}\mathcal N_A(x^\ast)v
		=cv+D\Tw\bigl((h'\circ x^\ast)v\bigr).
	\end{equation}
	The subspace $V_r$ is $\mathcal L$-invariant, and, in the block form
	induced by $L^2(I)=V_r\oplus V_r^\perp$,
	\begin{equation}\label{eq:block-A}
		\mathcal L=
		\begin{pmatrix}
			J & \ast\\[2pt]
			0 & cI
		\end{pmatrix}.
	\end{equation}

	\item[\textup{(ii)}]
	For class \textup{(B)}, suppose that
	$\tilde{g}(y,\cdot)\in C^1(\R)$ for a.e.\ $y\in I$ and that
$a(y):=
		\left.
		\partial_z\tilde{g}(y,z)
		\right|_{z=(\Tw x^\ast)(y)}
		\in L^\infty(I)$.
	Then, the vector field
	$\mathcal N_B(x)
		=cx+c_0\mathbf{1}+\tilde{g}(\cdot,\Tw x)$
	is Fr\'echet differentiable at $x^\ast$, and its derivative is the
	bounded operator
	\begin{equation}
		\mathcal L v
		:=\mathrm{D}\mathcal N_B(x^\ast)v
		=cv+a\cdot(\Tw v).
	\end{equation}
	The subspace $V_r^\perp$ is $\mathcal L$-invariant, and
	\begin{equation}\label{eq:block-B}
		\mathcal L=
		\begin{pmatrix}
			J & 0\\[2pt]
			\ast & cI
		\end{pmatrix}.
	\end{equation}
\end{enumerate}
In both cases, the $V_r\to V_r$ diagonal block is represented in the
basis $(\varphi_1,\ldots,\varphi_r)$ by the $r\times r$ matrix $J$.
This matrix is the Jacobian, evaluated at $X^\ast=\Mr x^\ast$, of the
reduced vector field in Theorem~\ref{thm:exact}:
\begin{equation}\label{eq:reduced-jacobian}
	J_{\mu\nu}=
	\begin{cases}
		c\delta_{\mu\nu}
		+D\lambda_\mu
		\displaystyle \int_0^1
\varphi_\mu(y)
h'\bigl(x^\ast(y)\bigr)
\varphi_\nu(y)\,\dd y
		& \textup{for class (A)},\\[3mm]
		c\delta_{\mu\nu}
		+\lambda_\nu
		\displaystyle \int_0^1
\varphi_\mu(y) a(y) \varphi_\nu(y)\,\dd y
		& \textup{for class (B)}.
	\end{cases}
\end{equation}
\end{lemma}

\begin{proof}

The boundedness of $\Tw:L^2(I)\to L^2(I)$, together with the
measurability and global Lipschitz assumptions in
Definition~\ref{def:reducible} and, for class \textup{(B)}, the
condition $\tilde g(\cdot,0)\in L^2(I)$, ensures that
$\mathcal N_A$ and $\mathcal N_B$ are well-defined maps from
$L^2(I)$ to $L^2(I)$.

We first verify that the operators stated in the lemma are the
Fr\'echet derivatives of the two vector fields. Afterward, we derive
their block forms with respect to $L^2(I)=V_r\oplus V_r^\perp$.

\smallskip
\noindent\emph{Class \textup{(A)}.}
The proposed derivative is
$\mathcal L v=cv+D\Tw\bigl((h'\circ x^\ast)v\bigr)$.
To verify this formula, we must show that the part of
$\mathcal N_A(x^\ast+v)-\mathcal N_A(x^\ast)$ not accounted for by
$\mathcal L v$ is $o(\norm{v}_{L^2})$. The only nonlinear contribution
to this difference comes from $h$. Therefore, we define its pointwise first-order
remainder by
\begin{equation}
R_v
=
h\circ(x^\ast+v)-h\circ x^\ast-(h'\circ x^\ast)v.
\end{equation}
Then,
\begin{equation}
\mathcal N_A(x^\ast+v)-\mathcal N_A(x^\ast)-\mathcal L v
=
D\Tw R_v.
\end{equation}
It therefore suffices to prove that
$\norm{\Tw R_v}_{L^2}
=
o(\norm{v}_{L^2})$.

Because smallness of $\norm{v}_{L^2}$ does not imply pointwise
smallness of $v$, we estimate $R_v$ separately on the region
where $\abs{v}\le\delta$, using uniform continuity of $h'$, and on the
region where $\abs{v}>\delta$, using the global Lipschitz bound on $h$
together with the small $L^2$ norm of $v$.

We start with the region in which $\abs{v}\le\delta$. Choose $\rho>0$ such
that the compact set
$\Omega_\rho
=
\{s\in\R:\operatorname{dist}(s,\Omega)\le\rho\}$
is contained in the neighborhood on which $h\in C^1$. For
$0<\delta\le\rho$, define the modulus of continuity
\begin{equation}
\omega(\delta)
=
\sup
\left\{
\abs{h'(s)-h'(t)}
: s,t\in\Omega_\rho,\ 
\abs{s-t}\le\delta
\right\}.
\end{equation}
Because $h'$ is continuous on the compact set $\Omega_\rho$, the
Heine--Cantor theorem guarantees that $h'$ is uniformly continuous on
$\Omega_\rho$. Therefore, we obtain
$\omega(\delta)\longrightarrow 0$ as $\delta\longrightarrow 0$.
On the set $\{\abs{v}\le\delta\}$, the fundamental theorem of calculus
gives
\begin{equation}
R_v
=
v\int_0^1
\left[
h'\bigl(x^\ast+t v\bigr)-h'(x^\ast)
\right]\dd t.
\end{equation}
Because $x^\ast\in S$ and therefore $x^\ast(y)\in\Omega$ for a.e.\
$y\in I$, all arguments of $h'$ in this integral belong to
$\Omega_\rho$. Thus, we obtain
\begin{equation}
\abs{R_v}
\le
\omega(\delta)\abs{v}
\label{eq:R_v-classA-1}
\end{equation}
on $\{\abs{v}\le\delta\}$.

We now consider the region in which $\abs{v} > \delta$.
Let $L_h$ be a global Lipschitz constant for $h$. 
On the set $\{\abs{v}>\delta\}$, global Lipschitz
continuity of $h$ and the bound
$\abs{h'(x^\ast)}\le L_h$ give
\begin{equation}
\abs{R_v}
\le
2L_h\abs{v}.
\label{eq:R_v-classA-2}
\end{equation}
Combining \eqref{eq:R_v-classA-1} and \eqref{eq:R_v-classA-2} yields
\begin{align}
\norm{R_v}_{L^1}
&\le
\omega(\delta)\norm{v}_{L^1}
+
2L_h
\int_{\{\abs{v}>\delta\}}
\abs{v}\,\dd y \notag\\
&\le
\omega(\delta)\norm{v}_{L^1}
+
2L_h
\int_{\{\abs{v}>\delta\}}
\frac{\abs{v}}{\delta} \cdot \abs{v}\,\dd y \notag\\
&\le
\omega(\delta)\norm{v}_{L^2}
+
\frac{2L_h}{\delta}\norm{v}_{L^2}^2.
\label{eq:R_v-classA-final-bound}
\end{align}

Because $0\le W\le1$ and $I$ has measure one, we obtain
\begin{equation}
\norm{\Tw u}_{L^2}
\le
\norm{u}_{L^1}.
\label{eq:Twu_L2-bounded-by-L1}
\end{equation}
By substituting \eqref{eq:R_v-classA-final-bound} as $u$ into \eqref{eq:Twu_L2-bounded-by-L1},
we obtain, for every fixed $\delta\in(0,\rho]$,
\begin{equation}
\limsup_{\norm{v}_{L^2}\to0}
\frac{\norm{\Tw R_v}_{L^2}}{\norm{v}_{L^2}}
\le
\omega(\delta).
\end{equation}
Letting $\delta\to0$ gives
\begin{equation}
\norm{\Tw R_v}_{L^2}
=
o(\norm{v}_{L^2}).
\end{equation}
Consequently,
\begin{equation}
\norm{
\mathcal N_A(x^\ast+v)
-
\mathcal N_A(x^\ast)
-
\mathcal L v
}_{L^2}
=
o(\norm{v}_{L^2}),
\end{equation}
which proves that
\begin{equation}
\mathrm{D}\mathcal N_A(x^\ast)v
=
cv+D\Tw\bigl((h'\circ x^\ast)v\bigr).
\end{equation}
This derivative is bounded because $h'\circ x^\ast\in L^\infty(I)$
(which follows because $x^\ast(y)\in\Omega$ for a.e.\
$y\in I$ and $h'$ is continuous, and therefore bounded, on the
compact set $\Omega$)
and $\Tw$ is bounded on $L^2(I)$.

\smallskip
\noindent\emph{Class \textup{(B)}.}
Let
$a(y)=\partial_z\tilde g\bigl(y,z^\ast(y)\bigr)$.
The proposed derivative is
$\mathcal L v=cv+a\cdot(\Tw v)$.
To verify this formula, set
$z^\ast=\Tw x^\ast$, $u_v=\Tw v$, 
and define the first-order remainder of $\tilde g$ in its second
argument by
\begin{equation}
Q_v
=
\tilde g(\cdot,z^\ast+u_v)
-
\tilde g(\cdot,z^\ast)
-
a\,u_v.
\end{equation}
Then,
\begin{equation}
\mathcal N_B(x^\ast+v)-\mathcal N_B(x^\ast)-\mathcal L v
=
Q_v.
\end{equation}
Thus, it remains to prove that
$\norm{Q_v}_{L^2}
=
o(\norm{v}_{L^2})$.

The graphon operator converts an $L^2$ perturbation into a uniformly
small perturbation of the second argument of $\tilde g$. Indeed, the
Cauchy--Schwarz inequality and $0\le W\le1$ give
\begin{equation}
\abs{u_v(y)}
\le
\left(
\int_0^1 W(y,y')^2\,\dd y'
\right)^{1/2}
\norm{v}_{L^2}
\le
\norm{v}_{L^2}
\end{equation}
for a.e.\ $y\in I$. Therefore, we obtain
\begin{equation}
\norm{u_v}_{L^\infty}
\le
\norm{v}_{L^2}.
\label{eq:u_v-v}
\end{equation}

Let $v^{(1)}$, $v^{(2)}$, $\ldots$ be any sequence of nonzero elements of $L^2(I)$ that converges to $0$ in $L^2(I)$. Then, $u_{v^{(m)}}(y)\to0$ as $m \to \infty$ for a.e.\ $y$. 
Whenever $u_{v^{(m)}}(y)\ne0$, we obtain
\begin{equation}
\frac{Q_{v^{(m)}}(y)}{\norm{v^{(m)}}_{L^2}}
=
\left[
\frac{
\tilde g\bigl(y,z^\ast(y)+u_{v^{(m)}}(y)\bigr)
-
\tilde g\bigl(y,z^\ast(y)\bigr)
}{
u_{v^{(m)}}(y)
}
-
a(y)
\right]
\frac{u_{v^{(m)}}(y)}{\norm{v^{(m)}}_{L^2}}
\longrightarrow 0
\label{eq:Q-pointwise-convergence}
\end{equation}
for a.e.\,$y\in I$
because
the expression in square brackets tends to zero 
and \eqref{eq:u_v-v} implies that
$\frac{\abs{u_{v^{(m)}}(y)}}{\norm{v^{(m)}}_{L^2}}
\le1$. The expression in the square brackets tends to zero because
$u_{v^{(m)}}(y)\to0$ and $\tilde g(y,\cdot)$ is differentiable at
$z^\ast(y)$ with derivative $a(y)$.

When $u_{v^{(m)}}(y)=0$, the definition of $Q_{v^{(m)}}$ gives
$Q_{v^{(m)}}(y)=0$, such that the same convergence conclusion holds.

Pointwise convergence (i.e., \eqref{eq:Q-pointwise-convergence}), which we have shown, alone does not imply convergence in $L^2$. So,
we next obtain an integrable bound, independent of $m$, to
apply the dominated convergence theorem.
Similar to the derivation of \eqref{eq:Q-pointwise-convergence}, the uniform Lipschitz condition on $\tilde g$ gives the bound
\begin{equation}
\frac{\abs{Q_{v^{(m)}}(y)}}{\norm{v^{(m)}}_{L^2}}
\le
L_{\tilde g}+\norm{a}_{L^\infty}
\qquad
\text{for a.e.\ }y\in I.
\end{equation}
The dominated convergence theorem applied to 
$\left(\frac{\left| Q_{v^{(m)}}(y) \right|}{\norm{v^{(m)}}_{L^2}} \right)^2$
yields
\begin{equation}
\frac{\norm{Q_{v^{(m)}}}_{L^2}}{\norm{v^{(m)}}_{L^2}}
\longrightarrow 0 \quad \text{ as } m \to\infty.
\label{eq:Q-L2-a}
\end{equation}
Because \eqref{eq:Q-L2-a} holds true for every sequence $v^{(m)}\to0$ in $L^2(I)$,
\begin{equation}
\norm{Q_v}_{L^2}
=
o(\norm{v}_{L^2}).
\end{equation}
Consequently,
\begin{equation}
\mathrm{D}\mathcal N_B(x^\ast)v
=
cv+a\cdot(\Tw v).
\end{equation}
This derivative is bounded because $a\in L^\infty(I)$ and $\Tw$ is
bounded on $L^2(I)$.

\smallskip
\noindent
Having established that the displayed operators are the
Fr\'echet derivatives, we now derive their block forms.

For class \textup{(A)}, if $v\in V_r$, then $cv\in V_r$ and
$\Tw\bigl((h'\circ x^\ast)v\bigr)
\in
\ran\Tw
=
V_r$.
Thus, $V_r$ is invariant under $\mathcal L$, which gives the zero
lower-left block in~\eqref{eq:block-A}. If $v\in V_r^\perp$, then
the first term of $\mathcal L v
=
cv+D\Tw\bigl((h'\circ x^\ast)v\bigr)$
belongs to $V_r^\perp$, and the second term belongs to
$V_r$. Therefore, the $V_r^\perp$-component of $\mathcal L v$ is
$cv$, which gives the lower-right block $cI$.

For $\mu,\nu\in\{1,\ldots,r\}$, self-adjointness of $\Tw$ and
$\Tw\varphi_\mu=\lambda_\mu\varphi_\mu$ give
\begin{equation}
\inner{\varphi_\mu}{\mathcal L\varphi_\nu}
=
c\delta_{\mu\nu}
+
D\lambda_\mu
\int_0^1
\varphi_\mu(y) h'(x^\ast(y))\varphi_\nu(y)\,\dd y.
\end{equation}
This is $\partial\dot X_\mu/\partial X_\nu$ computed
from~\eqref{eq:exactA}, so the $V_r\to V_r$ block is the class
\textup{(A)} matrix in~\eqref{eq:reduced-jacobian}.

For class \textup{(B)}, if $v\in V_r^\perp=\ker\Tw$, then
$\Tw v=0$, so $\mathcal L v=cv\in V_r^\perp$.
Thus, $V_r^\perp$ is invariant, $\mathcal L|_{V_r^\perp}=cI$, and the
upper-right block in~\eqref{eq:block-B} vanishes.

For $\mu,\nu\in\{1,\ldots,r\}$,
\begin{equation}
\inner{\varphi_\mu}{\mathcal L\varphi_\nu}
=
c\delta_{\mu\nu}
+
\lambda_\nu
\int_0^1
\varphi_\mu(y) a(y) \varphi_\nu(y)\,\dd y.
\end{equation}
This is $\partial\dot X_\mu/\partial X_\nu$ computed
from~\eqref{eq:exactB}, so the $V_r\to V_r$ block is the class
\textup{(B)} matrix in~\eqref{eq:reduced-jacobian}.
\end{proof}

The triangular linearization determines the spectrum and linear
stability. For local bifurcations, the correspondence follows from the
stronger exact nonlinear structures of the two classes, as shown below.

\begin{theorem}[The reduction preserves equilibria, stability, and bifurcations]\label{thm:bifurcation}
Let $W$ have finite rank $r$. Regard each vector field in Definition~\ref{def:reducible} (i.e., \eqref{eq:reducible-A} and \eqref{eq:reducible-B}) as a map
from $L^2(I)$ to $L^2(I)$.
For parts \textup{(ii)}--\textup{(iv)},
assume the hypotheses of Lemma~\ref{lem:triangular} at the equilibrium
under consideration. Whenever spectra are discussed, $\mathcal L$ and
$J$ are regarded, by complexification, as operators on
$L^2(I;\mathbb C)$ and $\mathbb C^r$, respectively, and $\spec$
denotes the corresponding complex spectrum.
\begin{enumerate}
\item[\textup{(i)}] \textup{(Equilibria.)} For class \textup{(A)}, the
observable map
\begin{equation}
x^\ast\longmapsto X^\ast=\Mr x^\ast
\end{equation}
is a bijection from the equilibria of the full graphon dynamics that
lie in $V_r$ onto the equilibria of the reduced
system~\eqref{eq:exactA}, with inverse
\begin{equation}
X^\ast\longmapsto x^\ast=\Mr^\ast X^\ast
=\sum_{\nu=1}^{r}X^\ast_\nu\varphi_\nu.
\end{equation}
If $c\neq0$, every equilibrium of the full class-\textup{(A)}
dynamics belongs to $V_r$, so this is a bijection between the entire
full and reduced equilibrium sets. If $c=0$, equilibria outside $V_r$
may exist, and the reduction does not in general represent them.

For class \textup{(B)} with $c\neq0$, the same observable map is a
bijection from all equilibria of the full graphon dynamics onto the
equilibria of~\eqref{eq:exactB}, with inverse
\begin{equation}\label{eq:equilibrium-bijection}
X^\ast\longmapsto x^\ast=-\frac1c\Bigl[c_0\mathbf{1}
+\tilde{g}\Bigl(\cdot,\ \textstyle\sum_{\nu=1}^{r}
\lambda_\nu X^\ast_\nu\varphi_\nu\Bigr)\Bigr].
\end{equation}
\item[\textup{(ii)}] \textup{(Spectrum.)}
$\spec(\mathcal L)=\spec(J)\cup\{c\}$.
Thus, the spectrum of the full linearization consists of the
eigenvalues of the reduced Jacobian together with the transverse
spectral point $c$, which may already belong to $\spec(J)$.
\item[\textup{(iii)}] \textup{(Stability.)} If $c<0$, then $x^\ast$ is linearly asymptotically stable for the graphon dynamics if and only if $X^\ast$ is linearly asymptotically stable for the reduced system. If $c\ge0$, no equilibrium is linearly asymptotically stable.
\item[\textup{(iv)}] \textup{(Local bifurcations.)}
Suppose that the model depends smoothly on a parameter $p$ through
$D$, $h$, or $\tilde g$, while $c<0$ is fixed. 
Let $X^\ast$ be an equilibrium of the reduced system at $p=p_0$. Let $x^\ast$ be its full lift from part \textup{(i)}. Let $J$ be the
Jacobian of the reduced vector field with respect to $X$, evaluated at
$(X^\ast,p_0)$. Finally, let $\mathcal L$ be the Fr\'echet derivative of
the full vector field with respect to $x$, evaluated at
$(x^\ast,p_0)$.
Then,
the full and reduced linearizations have exactly the
same spectral values with zero real part.

If the reduced system satisfies the smoothness, nondegeneracy, and
transversality conditions required for a local steady-state or Hopf
bifurcation at $(X^\ast,p_0)$, then its bifurcating equilibria or
periodic solutions lift uniquely to the full graphon dynamics.
Conversely, every nearby full equilibrium or periodic solution
projects to the corresponding reduced solution. Thus, the full and
reduced systems undergo the same local bifurcation at the same
parameter value and of the same type.

\end{enumerate}
\end{theorem}

\begin{proof}

(i) We first consider class \textup{(A)}. Let $P_r=\Mr^\ast\Mr$ be
the orthogonal projection onto $V_r$. If $x^\ast$ is a full
equilibrium, then
\begin{equation}
0=(I-P_r)\bigl[cx^\ast+D\Tw(h\circ x^\ast)\bigr]
=c(I-P_r)x^\ast,
\end{equation}
because $\Tw(h\circ x^\ast)\in\ran\Tw=V_r$. Hence, every full
equilibrium belongs to $V_r$ when $c\neq0$.

Now let $x^\ast\in V_r$ and write
$X^\ast=\Mr x^\ast$ and
$x^\ast=\Mr^\ast X^\ast
=\sum_{\nu=1}^rX_\nu^\ast\varphi_\nu$.
Taking the inner product of the full equilibrium equation with
$\varphi_\mu$ gives
\begin{equation}
0=cX_\mu^\ast
+D\lambda_\mu\int_0^1\varphi_\mu(y)
h\!\Bigl(\textstyle\sum_{\nu=1}^r
X_\nu^\ast\varphi_\nu(y)\Bigr)\dd y,
\label{eq:classA-equil-full}
\end{equation}
which is the equilibrium equation of~\eqref{eq:exactA}. Conversely, if
$X^\ast$ satisfies \eqref{eq:classA-equil-full} and
$x^\ast=\Mr^\ast X^\ast$, then the full vector field at $x^\ast$
belongs to $V_r$ and is orthogonal to every basis vector
$\varphi_\mu$. It therefore vanishes. This proves the stated
bijection for class \textup{(A)}.

For class \textup{(B)}, assume $c\neq0$. If $x^\ast$ is a full
equilibrium and $X^\ast=\Mr x^\ast$, then
$\Tw x^\ast
=
\sum_{\nu=1}^r\lambda_\nu X_\nu^\ast\varphi_\nu$.
By substituting this equality into
$cx^\ast+c_0\mathbf{1}+\tilde g(\cdot,\Tw x^\ast)=0$
and taking the inner product of it with each $\varphi_\mu$, we find that $X^\ast$ is an equilibrium
of~\eqref{eq:exactB}. Moreover, because $c\neq0$, the full equilibrium
is uniquely determined by $X^\ast$ as
\begin{equation}
x^\ast=-\frac1c\Bigl[c_0\mathbf{1}
+\tilde g\Bigl(\cdot,\textstyle\sum_{\nu=1}^r
\lambda_\nu X_\nu^\ast\varphi_\nu\Bigr)\Bigr].
\label{eq:x^*-from-X^*-classB}
\end{equation}

Conversely, let $X^\ast$ be an equilibrium of~\eqref{eq:exactB} and
define $x^\ast$ by \eqref{eq:x^*-from-X^*-classB}. The reduced equilibrium
equations imply that
$\inner{\varphi_\mu}{x^\ast}=X_\mu^\ast$,
$\mu \in \{ 1,\ldots,r \}$.
Hence, $\Mr x^\ast=X^\ast$ and
$\Tw x^\ast
=
\sum_{\nu=1}^r\lambda_\nu X_\nu^\ast\varphi_\nu$.
Substitution of the last equality into \eqref{eq:x^*-from-X^*-classB} shows that $x^\ast$
satisfies the full equilibrium equation. This proves the stated
bijection for class \textup{(B)}.

(ii) We work on the complexified spaces specified in the theorem,
using the same notation for the complexified operators and subspaces.
We identify $V_r$ with $\mathbb C^r$ through the basis
$(\varphi_1,\ldots,\varphi_r)$.

First, suppose that
$\lambda\notin\spec(J)\cup\{c\}$.
Then, the two diagonal blocks
$J-\lambda$ and $(c-\lambda)I$
of $\mathcal L-\lambda$ are boundedly invertible. A block triangular
operator whose diagonal blocks are boundedly invertible is itself
boundedly invertible. Therefore, 
$\lambda\notin\spec(\mathcal L)$,
which proves that
$\spec(\mathcal L)\subseteq\spec(J)\cup\{c\}$.

For the reverse inclusion, let $\lambda\in\spec(J)$. Because $J$ acts
on the finite-dimensional complex space $\mathbb C^r$, we find that $\lambda$ is
an eigenvalue of $J$.

For class \textup{(A)}, we choose a nonzero vector
$a=(a_1,\ldots,a_r)^\top\in\mathbb C^r$ such that
$Ja=\lambda a$,
and set
$w=\sum_{\nu=1}^r a_\nu\varphi_\nu\in V_r$.
Because $V_r$ is invariant under $\mathcal L$ and the matrix of
$\mathcal L|_{V_r}$ is $J$, we have
$\mathcal Lw
=
\sum_{\mu=1}^r(Ja)_\mu\varphi_\mu
=
\lambda w$.
Therefore, $\lambda\in\spec(\mathcal L)$.

For class \textup{(B)}, the upper-right block in
\eqref{eq:block-B} vanishes. Because $J-\lambda$ is not surjective,
choose $b\in V_r$ outside its range. If there were a
$v=v_\parallel+v_\perp$ satisfying
$(\mathcal L-\lambda)v=b$,
then taking the $V_r$ component would give
$(J-\lambda)v_\parallel=b$,
which contradicts that we have chosen $b$ to be outside the range of $J-\lambda$. Therefore, such a $v$ does not exist, and
$\mathcal L-\lambda$ is not surjective. Therefore,
$\lambda\in\spec(\mathcal L)$.

It remains to show that $c\in\spec(\mathcal L)$. For class
\textup{(A)}, the range of $\mathcal L-cI$ is contained in the proper
subspace $V_r$ of $L^2(I)$, so $\mathcal L-cI$ is not surjective. For class
\textup{(B)}, every nonzero $v\in V_r^\perp$ satisfies
$\mathcal Lv=cv$,
so $c$ is an eigenvalue of infinite multiplicity. Therefore,
$\spec(J)\cup\{c\}\subseteq\spec(\mathcal L)$,
which completes the proof of part \textup{(ii)}.

(iii) We write a perturbation from an equilibrium of the graphon dynamical system as $v(t)=v_\parallel(t)+v_\perp(t)$ with
$v_\parallel(t)\in V_r$ and $v_\perp(t)\in V_r^\perp$. For class \textup{(A)}, the block form
\eqref{eq:block-A} gives
$\dot v_\parallel=Jv_\parallel+Bv_\perp$ and $\dot v_\perp=cv_\perp$,
where $B:V_r^\perp\to V_r$ denotes the upper-right block. Therefore, we obtain
\begin{equation}
\begin{cases}
v_\parallel(t) &= e^{tJ}v_\parallel(0)+\int_0^t e^{(t-\tau)J}B e^{c\tau}v_\perp(0)\,\dd\tau,\\
v_\perp(t) &= e^{ct}v_\perp(0).
\end{cases}
\label{eq:stability-proof-A}
\end{equation}
For class \textup{(B)}, the block form \eqref{eq:block-B} gives
$\dot v_\parallel = Jv_\parallel$ and 
$\dot v_\perp =Bv_\parallel+cv_\perp$,
where now $B:V_r\to V_r^\perp$ denotes the lower-left block. Therefore, we obtain
\begin{equation}
\begin{cases}
v_\parallel(t) &= e^{tJ}v_\parallel(0),\\
v_\perp(t) &= e^{ct}v_\perp(0)+\int_0^t e^{c(t-\tau)}B e^{\tau J}v_\parallel(0)\,\dd\tau.
\end{cases}
\label{eq:stability-proof-B}
\end{equation}

Now assume that $c<0$ and that all eigenvalues of $J$ have negative real parts. Then, for both \eqref{eq:stability-proof-A} and \eqref{eq:stability-proof-B},
both the diagonal terms and the integral terms decay exponentially, such that the equilibrium is linearly asymptotically stable in $L^2(I)$. Conversely, if $J$ has an eigenvalue with nonnegative real part, then the dynamics of the $V_r$ component is not asymptotically stable due to the
$e^{tJ}$ factor. Therefore, the equilibrium of the full graphon dynamical system is not linearly asymptotically stable. 
If $c>0$, then $e^{ct}$ grows exponentially. If $c=0$, then the transverse component does not decay. So, if $c\ge0$, no equilibrium is linearly asymptotically stable.

(iv) By part \textup{(ii)}, we obtain
$\spec(\mathcal L)=\spec(J)\cup\{c\}$.
Because $c<0$, the full and reduced linearizations have exactly the same
spectral values with zero real part. We now use the exact nonlinear
structure of each class to relate their equilibria and periodic
solutions.

For class \textup{(A)}, let $P_r=\Mr^\ast\Mr$ and write
$u=P_r x$ and $z=(I-P_r)x$.
Because
$D\Tw\bigl(h\circ(u+z)\bigr)\in\ran\Tw=V_r$,
the full dynamics satisfies
\begin{equation}
\begin{cases}
\dot u &= cu+D\Tw\bigl(h\circ(u+z)\bigr),\\
\dot z &= cz.
\end{cases}
\end{equation}
At an equilibrium, $z=0$ because $c<0$. If $x(t)$ is $T$-periodic,
then $z(T)=e^{cT}z(0)=z(0)$.
Because $e^{cT}\ne1$, the equality $z(T) = z(0)$ implies $z(0)=0$ and hence $z(t)=0$ for all
$t$. Therefore, every equilibrium or periodic solution of the full
system lies in $V_r$. On this invariant subspace, the full dynamics
and the reduced dynamics are identical under the correspondence
$x=\Mr^\ast X$. Thus, their equilibria and periodic solutions
correspond uniquely.

For class \textup{(B)}, we write $X=\Mr x$ and $z=(I-P_r)x$.
Because
$\Tw x
=
\sum_{\nu=1}^r\lambda_\nu X_\nu\varphi_\nu$,
the full dynamics decomposes exactly as
\begin{equation}
\begin{cases}
\dot X &= F_r(X,p),\\
\dot z &= cz+H_p(X),
\label{eq:dyn-X-z-classB}
\end{cases}
\end{equation}
where $F_r$ is the reduced vector field and
$H_p(X)
=
(I-P_r)\left[
c_0\mathbf{1}
+
\tilde g_p\left(
\cdot,
\textstyle\sum_{\nu=1}^r
\lambda_\nu X_\nu\varphi_\nu
\right)
\right]$.
A reduced equilibrium $X^\ast$ has the unique full lift
$z^\ast=-\frac{1}{c}H_p(X^\ast)$.

Now let $X(t)$ be a $T$-periodic solution of the reduced system. Equation
\eqref{eq:dyn-X-z-classB} gives
\begin{equation}
z(t)
=
e^{ct}z(0)
+
\int_0^t e^{c(t-s)}H_p(X(s))\,\dd s.
\end{equation}
The condition $z(T)=z(0)$ uniquely determines
\begin{equation}
z(0)
=
\frac{1}{1-e^{cT}}
\int_0^T e^{c(T-s)}H_p(X(s))\,\dd s,
\end{equation}
because $e^{cT}\ne1$. Therefore, every reduced periodic solution has a
unique full periodic lift. Conversely, every full equilibrium or
periodic solution projects under $\Mr$ to an equilibrium or periodic
solution of the reduced system.

Therefore, under the bifurcation hypotheses stated in part
\textup{(iv)}, the nearby bifurcating equilibria or periodic solutions
of the reduced system correspond uniquely to those of the full
system. This correspondence preserves the parameter value and whether
the solutions are equilibria or nonconstant periodic solutions, which
proves the claimed bifurcation correspondence.

\end{proof}

\section{Approximation error for general graphon dynamics}\label{sec:error}

For nonlinear $f$ or $G$ and general (infinite-rank) $W$, the reduction~\eqref{eq:reduced} is approximate. In this section, we bound its error, isolating the roles of the state truncation error and, for input-coupling dynamics, the discarded singular values. We also show a Gr\"onwall estimate on the resulting trajectory error. Throughout, $x(\cdot,t)\in S$ denotes a fixed solution of the graphon dynamical system under consideration. We define the truncation error of the state by
\begin{equation}\label{eq:def-epsilon-n}
\varepsilon_n(x):=\norm{(I-\Pn)x}_{L^2(I)},
\end{equation}
where we remind that $\Pn$ denotes the projection onto $V_n$.

\subsection{Alignment accuracy}

A natural measure of the quality of the closure~\eqref{eq:reduced-op} is the alignment error in $\R^n$ between the exact observable dynamics~\eqref{eq:exact-obs} and the reduced field~\cite{thibeault2024low}. For the vector field $\mathcal N$ under consideration, we define
\begin{equation}\label{eq:align-def}
E_n(x):=\frac{1}{\sqrt n}\left\lVert \Mn\mathcal N(x)-F_n(\Mn x)\right\rVert_2
=\frac{1}{\sqrt n}\left\lVert \Mn\bigl(\mathcal N(x)-\mathcal N(\Pn x)\bigr)\right\rVert_2 ,
\end{equation}
where we used $F_n(\Mn x)=\Mn\mathcal N(\Mn^\ast\Mn x)=\Mn\mathcal N(\Pn x)$. We first treat input-coupling dynamics, for which the discarded modes enter only through $\abs{\lambda_{n+1}}$, exactly as in the finite-network theory in~\cite[Suppl.\ Thm.\ 72]{thibeault2024low}. This input-coupling form includes class \textup{(B)}, for which $f$ is affine, but the following estimate only requires $f$ to be Lipschitz.

\begin{theorem}[Alignment bound, input-coupling form]\label{thm:align}
Consider the input-coupling dynamics
\begin{equation}
\partial_t x=f(x)+\tilde{g}(\cdot,\Tw x),
\end{equation}
where $f$ is Lipschitz with constant $L_f$, $\tilde g(\cdot,0)\in L^2(I)$, and $\tilde{g}(y,\cdot)$ is Lipschitz with constant $L_{\tilde{g}}$, uniformly in $y$. For every state $x\in S$, we obtain
\begin{equation}\label{eq:align-bound}
\sqrt nE_n(x)
\le
\left(
L_f+L_{\tilde g}\abs{\lambda_{n+1}}
\right)
\varepsilon_n(x).
\end{equation}
\end{theorem}

\begin{proof}
The assumption $\tilde g(\cdot,0)\in L^2(I)$ and the Lipschitz condition in the second argument ensure that $\tilde g(\cdot,\Tw x)\in L^2(I)$ whenever $x\in L^2(I)$. 
Theefore, $\mathcal N(z)=f(z)+\tilde g(\cdot,\Tw z)$
belongs to $L^2(I)$ for every $z\in L^2(I)$, which leads to
$\mathcal N(x)-\mathcal N(\Pn x)\in L^2(I)$.

Bessel's inequality guarantees that $\norm{\Mn u}_2\le\norm{u}_{L^2(I)}$. Substitution of this inequality into \eqref{eq:align-def} yields $\sqrt nE_n(x)\le \norm{\mathcal N(x)-\mathcal N(\Pn x)}_{L^2(I)}$. Note that the right-hand side of this inequality is well-posed owing to
$\mathcal N(x)-\mathcal N(\Pn x)\in L^2(I)$.

To evaluate $\norm{\mathcal N(x)-\mathcal N(\Pn x)}_{L^2(I)}$, we split $\mathcal N=f+C$ with $C(x)=\tilde{g}(\cdot,\Tw x)$. For the self term, we obtain $\abs{f(x(y))-f(\Pn x(y))}\le L_f\abs{(I-\Pn)x(y)}$ pointwise, so $\norm{f(x)-f(\Pn x)}_{L^2}\le L_f\varepsilon_n(x)$.
For bounding the coupling term, we define
\begin{equation}\label{eq:def-delta}
\delta(y):=(\Tw x)(y)-(\Tw\Pn x)(y)=\bigl(\Tw(I-\Pn)x\bigr)(y).
\end{equation}
Using the Lipschitz continuity of $\tilde{g}(y,\cdot)$ and \eqref{eq:def-delta}, we obtain
\begin{equation}
\abs{C(x)(y)-C(\Pn x)(y)}=\abs{\tilde{g}\bigl(y,(\Tw x)(y)\bigr)-\tilde{g}\bigl(y,(\Tw\Pn x)(y)\bigr)}\le L_{\tilde{g}}\abs{\delta(y)} .
\end{equation}
Therefore, we obtain
\begin{align}
\norm{C(x)-C(\Pn x)}_{L^2}
&\le
L_{\tilde g}\norm{\Tw(I-\Pn)x}_{L^2} \notag \\
&\le
L_{\tilde g}\abs{\lambda_{n+1}}
\norm{(I-\Pn)x}_{L^2} \notag \\
&=
L_{\tilde g}\abs{\lambda_{n+1}}\varepsilon_n(x).
\end{align}
The second inequality follows from
Corollary~\ref{cor:restricted-norm}, because
$(I-\Pn)x\in V_n^\perp$. Adding this estimate to the self-term bound
gives~\eqref{eq:align-bound}.

\end{proof}

In~\eqref{eq:align-bound}, the coupling contribution is bounded by
$L_{\tilde g}\abs{\lambda_{n+1}}\varepsilon_n(x)$ and is therefore
small when both the discarded state component and
$\abs{\lambda_{n+1}}$ are small.
Here,
$\abs{\lambda_{n+1}}=\sigma_{n+1}$ is the $(n{+}1)$st singular value
of the graphon, in direct analogy with the finite-network estimate
in~\cite[Suppl.\ Thm.\ 72]{thibeault2024low}.

The general graphon dynamical system~\eqref{eq:graphon-ds}, in which $G$ depends on the state $x(y)$ at the point $y$ (not only on the input $\Tw x$), admits a coarser bound as follows.
\begin{proposition}[Alignment bound, general coupling]\label{prop:align-general}
For the general dynamics~\eqref{eq:graphon-ds} with Assumption~\ref{as:lip} and graphon degree $s(y)=(\Tw\mathbf{1})(y)$, let $s_\ast:=\operatorname{ess\,sup}_y s(y)\le1$. Then,
\begin{equation}\label{eq:align-general}
\sqrt nE_n(x)
\le
\bigl(L_f+L_G(s_\ast+\abs{\lambda_1})\bigr)\varepsilon_n(x),
\end{equation}
where $\lambda_1$ is the eigenvalue of $\Tw$ with the largest modulus.
\end{proposition}

\begin{proof}
As in the proof of Theorem~\ref{thm:align},
\begin{equation}
\sqrt n E_n(x)
\le
\norm{f(x)-f(\Pn x)}_{L^2}
+
\norm{C(x)-C(\Pn x)}_{L^2},
\end{equation}
where
\begin{equation}
C(x)(y)=\int_0^1 W(y,y')G(x(y),x(y'))\,\dd y' .
\end{equation}
The self term satisfies
\begin{equation}
\norm{f(x)-f(\Pn x)}_{L^2}
\le
L_f\norm{(I-\Pn)x}_{L^2}.
\end{equation}
Let $v=(I-\Pn)x$. For the coupling term, Assumption~\ref{as:lip} gives, for a.e.\ $y$,
\begin{align}
\abs{C(x)(y)-C(\Pn x)(y)}
&\le
L_G\int_0^1 W(y,y')\Bigl(\abs{v(y)}+\abs{v(y')}\Bigr)\,\dd y' \notag\\
&=
L_G\Bigl(s(y)\abs{v(y)}+(\Tw\abs{v})(y)\Bigr).
\end{align}
Taking $L^2$ norms yields
\begin{equation}
\norm{C(x)-C(\Pn x)}_{L^2}
\le
L_G \left( s_\ast\norm{v}_{L^2}
+
\norm{\Tw\abs{v}}_{L^2} \right).
\end{equation}
Because $\abs{v}$ has the same $L^2$ norm as $v$, we obtain
\begin{equation}
\norm{\Tw\abs{v}}_{L^2}
\le
\norm{\Tw}_{\mathrm{op}}\norm{\abs{v}}_{L^2}
=
\norm{\Tw}_{\mathrm{op}}\norm{v}_{L^2}.
\end{equation}
By Proposition~\ref{prop:svd}, $\norm{\Tw}_{\mathrm{op}}=\sigma_1=\abs{\lambda_1}$. Combining the estimates gives~\eqref{eq:align-general}.
\end{proof}

\subsection{Trajectory accuracy}

The alignment error measures the instantaneous mismatch of the vector fields. We now bound the error between the true observables
\begin{equation}
X^{\mathrm{true}}_\mu(t) = \inner{\varphi_\mu}{x(\cdot,t)}
\label{eq:X^true_mu}
\end{equation}
and the solution $X(t)$ of the reduced dynamical system~\eqref{eq:reduced-op}. Because $S$ consists of the states taking values in the bounded interval $\Omega$, every $u\in S$ satisfies $\norm{u}_{L^2(I)}\le B_\Omega$. ($B_\Omega$ is defined in
\eqref{eq:B_Omega-def}.) Therefore, Bessel's inequality gives $\norm{\Mn u}_2\le B_\Omega$. Thus, the true observables lie in the closed ball
$\mathcal{B}_n:=\bigl\{\xi\in\R^n:\ \norm{\xi}_2\le B_\Omega\bigr\}$.
The trajectory of the reduced dynamical system need not remain in this ball a priori. Therefore, in the estimate below, we use a Lipschitz constant for $F_n$ on a set containing both trajectories. 
With the Lipschitz extensions of $f$ and $G$ fixed at the beginning of Section~\ref{sec:reduction}, one can take such a constant globally on $\R^n$.

\begin{theorem}[Trajectory error]\label{thm:accuracy}
Let $x(\cdot,t)\in S$ be the solution of the graphon dynamical system on $[0,T]$. Let $X(t)\in\R^n$ be the solution of the reduced dynamical system $\dot X=F_n(X)$ (see~\eqref{eq:reduced-op}). Assume that $F_n$ is Lipschitz on $\R^n$ with Lipschitz constant $\kappa_n$. Define the residual
\begin{equation}
\rho_n(t):=\sqrt nE_n\bigl(x(\cdot,t)\bigr).
\end{equation}
Then, for all $t\in[0,T]$, we obtain
\begin{equation}\label{eq:gronwall}
\bigl\lVert X^{\mathrm{true}}(t)-X(t)\bigr\rVert_2
\le
e^{\kappa_n t}\bigl\lVert X^{\mathrm{true}}(0)-X(0)\bigr\rVert_2
+
\int_0^t e^{\kappa_n(t-\tau)}\rho_n(\tau)\,\dd\tau .
\end{equation}
In particular, if $X(0)=X^{\mathrm{true}}(0)=\Mn x(\cdot,0)$ and $\rho_n(\tau)\le\bar\rho_n$ on $[0,T]$, then
\begin{equation}
\bigl\lVert X^{\mathrm{true}}(t)-X(t)\bigr\rVert_2
\le
\bar\rho_n\frac{e^{\kappa_n t}-1}{\kappa_n},
\label{eq:trajectory-accuracy-special-case}
\end{equation}
with the right-hand side interpreted as $\bar\rho_n t$ when $\kappa_n=0$.
\end{theorem}

\begin{proof}
Let
\begin{equation}
e(t):=X^{\mathrm{true}}(t)-X(t).
\end{equation}
Since $\Mn$ is linear and time-independent, using \eqref{eq:X^true_mu}, we obtain
\begin{align}
\dot X^{\mathrm{true}}
&=
\left(\inner{\varphi_1}{\dot{x}(\cdot,t)},\ldots,\inner{\varphi_n}{\dot{x}(\cdot,t)}\right) \notag\\
&=
\left(\inner{\varphi_1}{\mathcal N(x)},\ldots,\inner{\varphi_n}{\mathcal N(x)}\right) \notag\\
&=
\Mn\mathcal N(x).
\end{align}
Therefore,
\begin{align}
\dot e
&=
\Mn\mathcal N(x)-F_n(X) \notag\\
&=
\bigl(\Mn\mathcal N(x)-F_n(X^{\mathrm{true}})\bigr)
+
\bigl(F_n(X^{\mathrm{true}})-F_n(X)\bigr).
\label{eq:gronwall-2}
\end{align}
The first bracketed term in~\eqref{eq:gronwall-2} is the alignment residual:
\begin{equation}
\Mn\mathcal N(x)-F_n(X^{\mathrm{true}})
=
\Mn\mathcal N(x)-F_n(\Mn x)
=
\Mn\mathcal N(x)-\Mn\mathcal N(\Pn x),
\end{equation}
whose Euclidean norm is $\rho_n(t)=\sqrt nE_n(x(\cdot,t))$ by~\eqref{eq:align-def}. Integrating~\eqref{eq:gronwall-2} from $0$ to $t$ and taking norms gives
\begin{align}
\norm{e(t)}_2
&\le
\norm{e(0)}_2
+
\int_0^t \rho_n(\tau)\,\dd\tau
+
\int_0^t \norm{F_n(X^{\mathrm{true}}(\tau))-F_n(X(\tau))}_2\,\dd\tau \notag\\
&\le
\norm{e(0)}_2
+
\int_0^t \rho_n(\tau)\,\dd\tau
+
\kappa_n\int_0^t \norm{e(\tau)}_2\,\dd\tau .
\label{eq:gronwall-integral}
\end{align}
Gronwall's inequality applied to~\eqref{eq:gronwall-integral} yields~\eqref{eq:gronwall}.
%
%
The bound \eqref{eq:trajectory-accuracy-special-case} follows immediately.

\end{proof}

\begin{corollary}[Reconstructed-state error]
\label{cor:reconstructed-state-error}
Under the hypotheses of Theorem~\ref{thm:accuracy}, define the
reconstructed state by
$\widehat x(\cdot,t):=\Mn^\ast X(t)$.
Then, for every $t\in[0,T]$,
\begin{equation}\label{eq:reconstructed-state-error}
\norm{x(\cdot,t)-\widehat x(\cdot,t)}_{L^2}^2
=
\varepsilon_n\bigl(x(\cdot,t)\bigr)^2
+
\norm{X^{\mathrm{true}}(t)-X(t)}_2^2.
\end{equation}
In particular, if $X(0)=X^{\mathrm{true}}(0)$, then
\begin{equation}\label{eq:reconstructed-state-bound}
\norm{x(\cdot,t)-\widehat x(\cdot,t)}_{L^2}
\le
\left[
\varepsilon_n\bigl(x(\cdot,t)\bigr)^2
+
\left(
\int_0^t
e^{\kappa_n(t-\tau)}\rho_n(\tau)\,\dd\tau
\right)^2
\right]^{1/2}.
\end{equation}
\end{corollary}

\begin{proof}
Because $\Pn=\Mn^\ast\Mn$ and
$X^{\mathrm{true}}(t)=\Mn x(\cdot,t)$, we obtain
\begin{align}
x(\cdot,t)-\widehat x(\cdot,t)
&= (I-\Pn)x(\cdot,t) + \Pn x(\cdot,t) - \Mn^\ast X(t)\notag\\
&=
(I-\Pn)x(\cdot,t) +
\Mn^\ast\bigl[X^{\mathrm{true}}(t)-X(t)\bigr].
\label{eq:reconstruction-state-error-proof-1}
\end{align}
The first term on the right-hand side of \eqref{eq:reconstruction-state-error-proof-1} belongs to $V_n^\perp$, whereas the second term belongs to
$V_n$, so the two terms are orthogonal. Moreover, $\Mn^\ast$ is an
isometry from $\R^n$ to $V_n$. The Pythagorean theorem therefore
gives~\eqref{eq:reconstructed-state-error}. The bound
\eqref{eq:reconstructed-state-bound} follows from
Theorem~\ref{thm:accuracy}.
\end{proof}

\begin{remark}
Under Assumption~\ref{as:lip} and the Lipschitz extensions of the nonlinearities fixed at the beginning of Section~\ref{sec:reduction}, the required global Lipschitz property of $F_n$ holds. For the general graphon dynamical system~\eqref{eq:graphon-ds}, one may take
\begin{equation}
\kappa_n\le L_f+L_G(s_\ast+\abs{\lambda_1}),
\label{eq:kappa_n-1}
\end{equation}
independently of $n$. For input-coupling dynamics $\partial_t x=f(x)+\tilde g(\cdot,\Tw x)$, one may take
\begin{equation}
\kappa_n\le L_f+L_{\tilde g}\abs{\lambda_1}.
\label{eq:kappa_n-2}
\end{equation}

We obtain \eqref{eq:kappa_n-1} and \eqref{eq:kappa_n-2} as follows. Because $\Mn$ and $\Mn^\ast$ have operator norm $1$, we obtain
\begin{equation}
\norm{F_n(\xi)-F_n(\eta)}_2
\le
\norm{\mathcal N(\Mn^\ast\xi)-\mathcal N(\Mn^\ast\eta)}_{L^2(I)}.
\end{equation}
For the graphon dynamical system~\eqref{eq:graphon-ds}, the estimate used in the proof of Proposition~\ref{prop:align-general} gives, for arbitrary $u,v\in L^2(I)$,
\begin{align}
\norm{\mathcal N(u)-\mathcal N(v)}_{L^2}
&\le
\left[ L_f+L_G(s_\ast+\norm{\Tw}_{\mathrm{op}})\right] \norm{u-v}_{L^2} \notag\\
&=
\left[ L_f+L_G(s_\ast+\abs{\lambda_1})\right] \norm{u-v}_{L^2}.
\end{align}
This gives \eqref{eq:kappa_n-1}. For input-coupling dynamics, we obtain
\begin{equation}
\norm{\mathcal N(u)-\mathcal N(v)}_{L^2}
\le
L_f\norm{u-v}_{L^2}
+
L_{\tilde g}\norm{\Tw(u-v)}_{L^2}
\le
\bigl(L_f+L_{\tilde g}\abs{\lambda_1}\bigr)\norm{u-v}_{L^2},
\end{equation}
which gives the second bound.
\end{remark}

Combining Theorem~\ref{thm:accuracy} with the alignment bound~\eqref{eq:align-bound} gives, for input-coupling dynamics started at $X(0)=\Mn x(0)$,
\begin{equation}\label{eq:combined}
\bigl\lVert X^{\mathrm{true}}(t)-X(t)\bigr\rVert_2
\le
\int_0^t e^{\kappa_n(t-\tau)}
\left(
L_f+L_{\tilde g}\abs{\lambda_{n+1}}
\right)
\varepsilon_n\bigl(x(\tau)\bigr)
\dd\tau ,
\end{equation}
which bounds the trajectory error in terms
of (a) the state truncation error $\varepsilon_n$ and (b) the tail of
the graphon spectrum through $\abs{\lambda_{n+1}}$. This estimate is
generally not sharp for linear dynamics. In the linear system of
Corollary~\ref{cor:linear-exact}, both terms on the right-hand side
may be nonzero, whereas the left-hand side is exactly zero when the
initial observables are matched.

\section{Numerical illustration of the spectral reduction}\label{sec:examples}

In this section, we illustrate the theory numerically. We use two graphon dynamical systems and three benchmark graphons, one of finite rank and two of infinite rank.

\subsection{Two graphon neural-field dynamical systems}\label{sec:nf-examples}

We use two neural-field models of interacting neuronal populations. They differ in whether the sigmoidal nonlinearity acts on the presynaptic activity before integration over the network, or on the total synaptic input after integration.

For the first dynamical system, the sigmoid acts on the presynaptic activity
$x(y')$ before integration~\cite{amari1977dynamics}. The corresponding graphon dynamical system is given by
\begin{equation}\label{eq:classA-example}
\partial_t x(y,t)
=
-x(y,t)
+ \tilde f(y)
+
D\int_0^1
W(y,y')\,
\frac{1}{1+ e^{\mu-\delta x(y',t)}}
\,\dd y' ,
\end{equation}
where $\tilde f\in L^2(I)$ is a fixed external input, $D$ is the coupling strength, $\mu$ and $\delta$ are parameters, and $\delta > 0$. Previous literature on dimension reduction~\cite{laurence2019spectral, bramburger2024persistence, Eldo2025arxiv} omitted $\tilde{f}$.
When
$\tilde f\equiv0$, equation~\eqref{eq:classA-example} belongs to class
\textup{(A)} with $c=-1$ and $h(u)=1/(1+e^{\mu-\delta u})$.
When $\tilde f\ne0$, it is an additively forced extension of class
\textup{(A)}. For simplicity, we refer to \eqref{eq:classA-example} as the class
\textup{(A)} model. On a graphon of finite rank $r$, if $\tilde f\in V_r$ and the initial
condition lies in $V_r$, then the forced vector field maps $V_r$ into
itself and the reduction at dimension $r$ is exact. This includes the
unforced case $\tilde f\equiv0$. On an infinite-rank graphon, the
reduction is approximate and its trajectory error is controlled by
Theorem~\ref{thm:accuracy}.

%
%

As a class-\textup{(B)} dynamics, we consider the Wilson--Cowan convention~\cite{wilson1973mathematical}, in which the sigmoid function acts on the total synaptic input. It is the graphon analogue of an additive recurrent neural network whose nonlinearity acts on the total weighted input, the form whose low-rank reduction was studied on finite networks in~\cite{thibeault2024low}. We adapt this additive form to the graphon setting as
\begin{equation}\label{eq:classB-example}
\partial_t x(y,t)
=
-x(y,t)
+
\frac{1}{1+\exp\!\left[\mu-\delta\left(\tilde f(y)+D\displaystyle\int_0^1 W(y,y')x(y',t)\,\dd y'\right)\right]},
\end{equation}
where $\tilde f\in L^2(I)$ is again the external input. Equation~\eqref{eq:classB-example} takes the form $\partial_t x(y,t)=cx(y,t)+c_0+\tilde g\bigl(y,(\Tw x)(y,t)\bigr)$ with $c=-1$, $c_0=0$, and $\tilde g(y,z)=1/\{1+\exp[\mu-\delta(\tilde f(y)+Dz)]\}$, so it is of class \textup{(B)} even when $\tilde f\neq0$. The classical Wilson--Cowan equation additionally multiplies this input nonlinearity by a refractory factor that is approximately equal to $1-x(y,t)$~\cite{wilson1973mathematical}.
%
%
That factor makes the coupling depend on the local activity $x(y,t)$ in addition to the aggregated input $(\Tw x)(y,t)$, which takes the model outside class \textup{(B)}. We therefore omit it and use~\eqref{eq:classB-example}, which belongs to class \textup{(B)} and is a standard firing-rate model in its own right. The coupling function $\tilde g(y,\cdot)$ is Lipschitz uniformly in $y$ with constant $\abs{D}\delta/4$,
%
%
the factor $\delta/4$ being the maximum slope of the sigmoid. Because the sigmoid is bounded and the self-dynamics $-x$ is inward-pointing, a bounded invariant interval exists, so the standing assumptions hold. On a finite-rank graphon of rank $r$, Theorem~\ref{thm:exact}(ii) closes the reduced system exactly at dimension $r$ for every initial condition. On an infinite-rank graphon, the discarded modes enter the input-coupling error only through $\abs{\lambda_{n+1}}$ by Theorem~\ref{thm:align}.

In both models, we set $\mu=3$ and $\delta=1$~\cite{laurence2019spectral} and use the coupling strength $D$ as the control parameter.

\subsection{Graphons and their spectra}\label{sec:graphons-spectra}

We use the following three graphons.
First, the two-block modular graphon has blocks $B_1=[0,\gamma)$ and $B_2=[\gamma,1]$ and takes the value $p_{\mathrm{in}}$ within a block and $p_{\mathrm{out}}$ between blocks. 
It is constant on the $2\times2$ block partition, so $\Tw$ maps every function to a block-constant function and hence has rank at most two~\cite{klimm2022modularity}.
Its eigenfunctions are constant on each block, and the nonzero eigenvalues are those of the $2\times2$ matrix
\begin{equation}
M=\begin{pmatrix}
p_{\mathrm{in}}\gamma & p_{\mathrm{out}}(1-\gamma)\\
p_{\mathrm{out}}\gamma & p_{\mathrm{in}}(1-\gamma)
\end{pmatrix},
\end{equation}
which acts on block-constant functions. The two eigenvalues are
$\lambda_{1,2} =
\left[
p_{\mathrm{in}}
\pm
\sqrt{
p_{\mathrm{in}}^2
-
4\gamma(1-\gamma)(p_{\mathrm{in}}^2-p_{\mathrm{out}}^2)
}
\right]\big/ 2$,
and they are generically both nonzero, so $\rank(W)=2$. We set
$\gamma=0.3,\qquad p_{\mathrm{in}}=0.2$, and $p_{\mathrm{out}}=0.01$,
for which $\lambda_1=0.140$ and $\lambda_2=0.060$. 
Therefore, the class-\textup{(B)} dynamics has an exact
two-dimensional reduction on this graphon for every initial
condition. The class-\textup{(A)} dynamics also reduces exactly to
two dimensions when the initial state and any additive forcing belong
to the two-dimensional block-constant subspace $V_2$.
The $L^2$-normalized eigenfunctions $\varphi_1,\varphi_2$ are block-constant. We find
$\varphi_1=(0.104,1.193)$ and
$\varphi_2=(1.823,-0.068)$, where $\varphi_k=(\varphi_{k,1},\varphi_{k,2})$ for the values of $\varphi_k$ on $(B_1,B_2)$.

Second, the small-world graphon~\cite{medvedev2014nonlinear, kuehn2019power} 
is a graphon variant of the Watts-Strogatz small-world network model
\cite{watts1998collective}. It is defined by
$W=(1-p)W_q+p(1-W_q)$, where $W_q$ is the ring graphon of radius $q$, defined by $W_q(y,y')=1$ if the circular distance between $y$ and $y'$ satisfies $\min\{\abs{y-y'},1-\abs{y-y'}\}\le q$, and $W_q(y,y')=0$ otherwise. Because $W$ is translation invariant on the circle, it is diagonalized by the Fourier basis. The eigenvalue associated with the constant eigenfunction $\mathbf{1}$ is
$p+(1-2p)2q$. In addition, for each $k\ge1$, the eigenvalue
$(1-2p)\sin(2\pi k q) / \pi k$
has multiplicity two, with associated eigenfunctions
$\sqrt2\cos(2\pi ky)$ and $\sqrt2\sin(2\pi ky)$. We set
$p=0.1$ and $q=1/10$. For these parameter values, we obtain $\rank(W)=\infty$, with eigenvalues decaying algebraically as $O(1/k)$. 
Ordering the eigenvalues by decreasing modulus and indexing them from $\lambda_1$, the leading values are
$\lambda_1=0.260$ (associated with eigenfunction $\mathbf{1}$), $\lambda_2=\lambda_3=0.150$, $\lambda_4=\lambda_5=0.121$, $\lambda_6=\lambda_7=0.081$, and $\lambda_8=\lambda_9=0.037$.

\begin{figure}[htbp]
\centering
\includegraphics[width=0.6\textwidth]{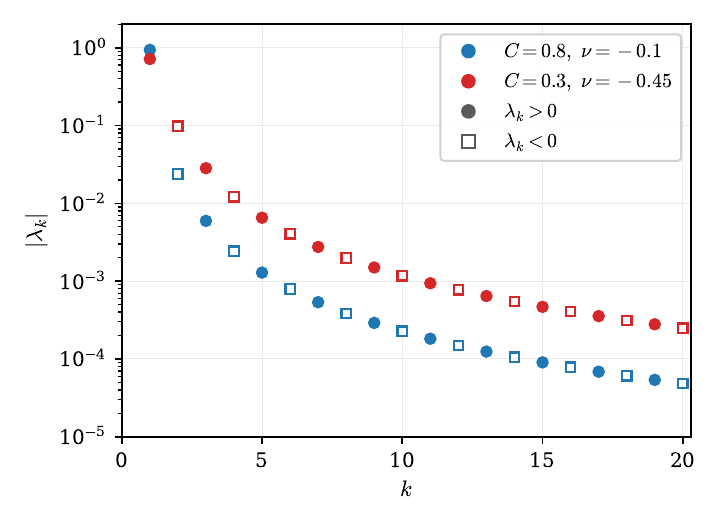}
\caption{Eigenvalue spectrum of the capped product-power graphon $W(y,y')=\min\{C(yy')^\nu,1\}$. We plot $\abs{\lambda_k}$ in descending order of modulus for two parameter sets. Filled circles mark positive eigenvalues. Open squares mark negative eigenvalues. Blue and red distinguish the two parameter sets. We discretized $[0,1]$ into $N=2000$ equally sized bins to compute the spectrum.}
\label{fig:pl-spectrum}
\end{figure}

Third, we consider the capped product-power graphon
given by $W(y,y')= 1$ if $y=0$ or $y' = 0$, and $W(y,y') = \min\{C(yy')^\nu,1\}$ otherwise,
where $0<C<1$ and $-0.5<\nu<0$.
This or similar graphon were also used
in~\cite{van2018sparse, medvedev2020kuramoto, Eldo2025arxiv}
and is related to low-rank network models with power-law degree distributions \cite{thibeault2024low}.
This graphon realizes a bounded heterogeneous degree distribution.
We compute the spectrum of the capped product-power graphon numerically.
Figure~\ref{fig:pl-spectrum} shows $\abs{\lambda_k}$ in descending order
of modulus for two parameter sets, $(C,\nu)=(0.8,-0.1)$ and
$(C,\nu)=(0.3,-0.45)$. In both cases, there is a single dominant
eigenvalue followed by substantially smaller eigenvalues of both signs.
For these two parameter sets, the ratios $\abs{\lambda_2/\lambda_1}$
are $0.025$ and $0.137$, respectively; thus, the spectrum is more
strongly dominated by its leading mode for $(C,\nu)=(0.8,-0.1)$.
We use this parameter set in the remainder of this section.

\subsection{Numerical methods}\label{sec:num-methods}

We discretize $I=[0,1]$ into $N$ bins with midpoints $y_i=(i-0.5)/N$, $i\in\{1,\ldots,N\}$. Approximating the coupling integral by the midpoint quadrature rule turns the graphon operator into the matrix with entries $W(y_i,y_j)/N$. We use the eigenvalues and eigenvectors of this matrix as discrete approximations of the eigenvalues and eigenfunctions of $\Tw$. We take $N=2000$ when computing spectra and $N=1000$ when simulating the dynamics. We integrate all the graphon dynamical systems up to $t=500$, which is far beyond the observed relaxation time in every case. We then use the final state as the numerical equilibrium.

The full $N$-dimensional system and the reduced $n$-dimensional system are integrated from consistent initial conditions, with the reduced initial condition obtained by projecting the full initial condition onto the first $n$ eigenfunctions. For the modular graphon, following~\cite{Eldo2025arxiv}, we use two spatially homogeneous initial conditions: a lower state $x(y,0)\equiv x_{\mathrm{low}}$ and an upper state $x(y,0)\equiv x_{\mathrm{high}}$, lying in the basins of the low- and high-activity branches, respectively. We take $x_{\mathrm{low}}=0$ for both models, $x_{\mathrm{high}}=8$ for the class-\textup{(A)} model, and $x_{\mathrm{high}}=1$ for the class-\textup{(B)} model; the interval $[0,1]$ is positively invariant for the
class-\textup{(B)} model.

We set the external input $\tilde f$ as follows. On the two-block modular graphon, we set $\tilde f\equiv0$ for both dynamical systems, so that the class-\textup{(A)} and class-\textup{(B)} dynamics are exactly the reducible forms of Definition~\ref{def:reducible} and the finite-rank reduction is exact. On the small-world and capped product-power graphons, we instead use the fixed heterogeneous input
$\tilde f(y)=1.2\left[\sin(2\pi y)+0.7\cos(2\pi y)+0.5\cos(4\pi y)\right]$
in both~\eqref{eq:classA-example} and~\eqref{eq:classB-example}, starting from $x(y,0)\equiv0$. Without such heterogeneity, the small-world equilibria are homogeneous because the small-world graphon has constant degree, and the reduction would be trivially exact at $n=1$. With $\tilde f\neq0$, the class-\textup{(B)} model remains of class \textup{(B)}, because the input enters through $\tilde g$. In contrast, a general additive heterogeneous input takes it outside class \textup{(A)}. This does not affect the alignment or trajectory-error estimates of Section~\ref{sec:error}, because a fixed additive input cancels in the vector-field differences controlled by those estimates. Its role here is to shift the equilibrium and excite higher modes against which the reduction is tested.

We quantify the approximation error as follows. For a reduction dimension $n$ and coupling strength $D$, we integrate the full graphon dynamical system to its equilibrium $x_D^\ast$ and the $n$-dimensional reduced system to its equilibrium observable vector $X_D^\ast=(X^\ast_{D,1},\ldots,X^\ast_{D,n})$. We reconstruct
$\widehat x_D^\ast=\sum_{k=1}^{n}X^\ast_{D,k}\varphi_k$
and compute the relative $L^2$ error
$e(n,D)=\norm{\widehat x_D^\ast-x_D^\ast}_{L^2} / \norm{x_D^\ast}_{L^2}$.
By Corollary~\ref{cor:reconstructed-state-error}, the squared
numerator is the sum of the squared state truncation error and the
squared error in the retained equilibrium observables.
Finally, we aggregate this error over a grid $\mathcal D$ of coupling strengths by
\begin{equation}\label{eq:error-measure}
\overline{E}(n)=\Bigl[\tfrac{1}{\abs{\mathcal D}}\textstyle\sum_{D\in\mathcal D}e(n,D)^2\Bigr]^{1/2}.
\end{equation}

\subsection{Numerical results}\label{sec:num-results}

\begin{figure}[htbp]
\centering
\includegraphics[width=0.82\textwidth]{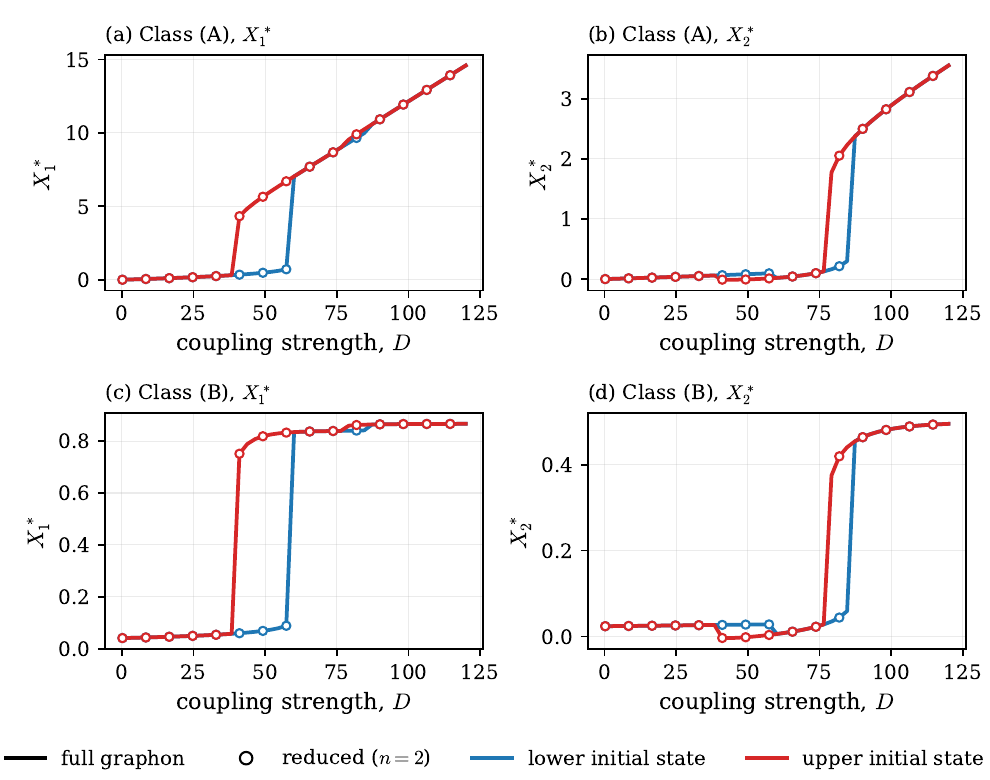}
\caption{Exact two-dimensional reduction on a modular graphon with two blocks. Panels (a) and (b) show the equilibrium observables $X_1^\ast$ and $X_2^\ast$, respectively, for the class-\textup{(A)} dynamics~\eqref{eq:classA-example}. Panels (c) and (d) show $X_1^\ast$ and $X_2^\ast$, respectively, for the class-\textup{(B)} dynamics~\eqref{eq:classB-example}. Solid lines show the observables computed from the full graphon dynamics. Open circles show the equilibria of the two-dimensional spectral reduction.}
\label{fig:modular-num}
\end{figure}

Figure~\ref{fig:modular-num} compares, as a function of the coupling
strength $D$, the two equilibrium observables $X_1^\ast$ and
$X_2^\ast$ of the full graphon dynamics (solid lines) with those of
the two-dimensional reduction (open circles), starting from lower and
upper initial states, for both graphon dynamical systems. For the
stable equilibria reached from these initial conditions, the full and
reduced observables agree to numerical precision throughout the
sampled values of $D$, including the associated hysteretic jumps, in
agreement with Theorems~\ref{thm:exact} and~\ref{thm:bifurcation}.
The observed hysteresis is consistent with saddle-node bifurcations.
Because the two blocks have different sizes, they switch at different
values of $D$, producing a two-stage cascade.

\begin{figure}[htbp]
\centering
\includegraphics[width=0.78\textwidth]{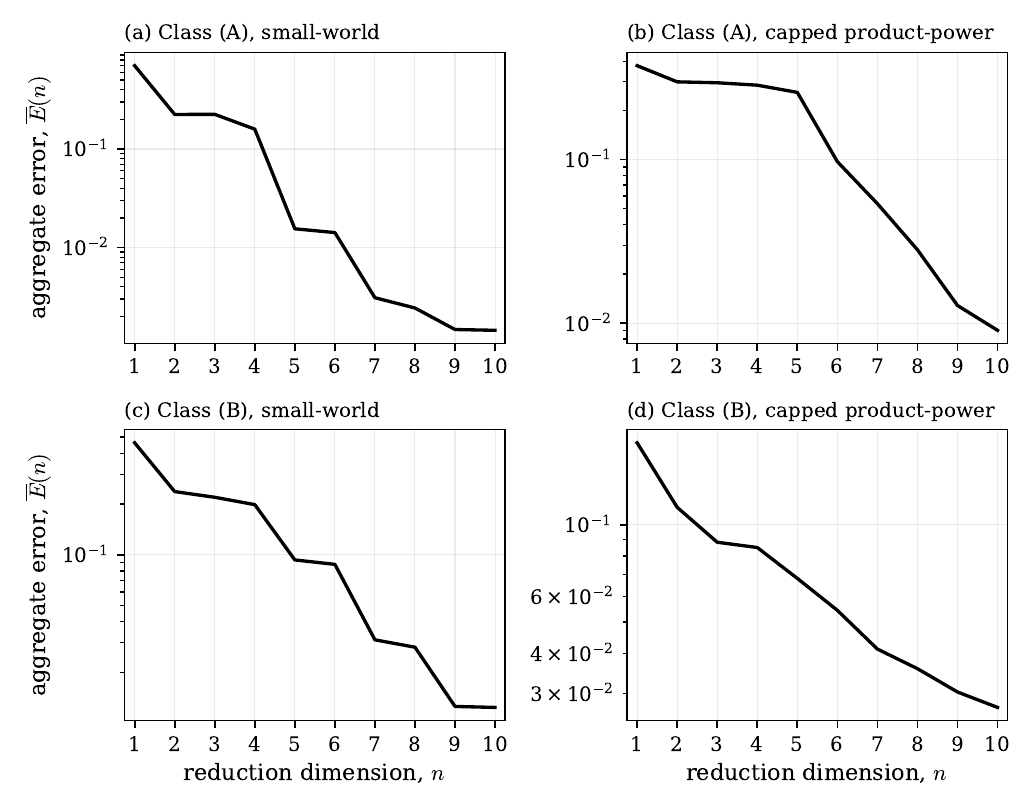}
\caption{Aggregate approximation error, \eqref{eq:error-measure}, for the reduced dynamical systems on the small-world and capped product-power graphons. (a) Class-\textup{(A)} dynamics on the small-world graphon. (b) Class-\textup{(A)} dynamics on the capped product-power graphon. (c) Class-\textup{(B)} dynamics on the small-world graphon. (d) Class-\textup{(B)} dynamics on the capped product-power graphon.}
\label{fig:error-num}
\end{figure}

Figure~\ref{fig:error-num} shows the aggregate
error~\eqref{eq:error-measure} as a function of the reduction
dimension $n$ for the small-world and capped product-power graphons
and for both dynamical systems. Overall, the error decreases as eigenfunctions
are added. For the
small-world graphon, whose spectrum decays algebraically as $O(1/k)$,
the error generally decreases in a staircase pattern as successive
Fourier modes are included
(Fig.~\ref{fig:error-num}(a) and (c)). For the capped product-power
graphon, whose spectrum has a single dominant mode followed by
substantially smaller modes, the error decreases quickly
(Fig.~\ref{fig:error-num}(b) and (d)).

\begin{figure}[htbp]
\centering
\includegraphics[width=0.82\textwidth]{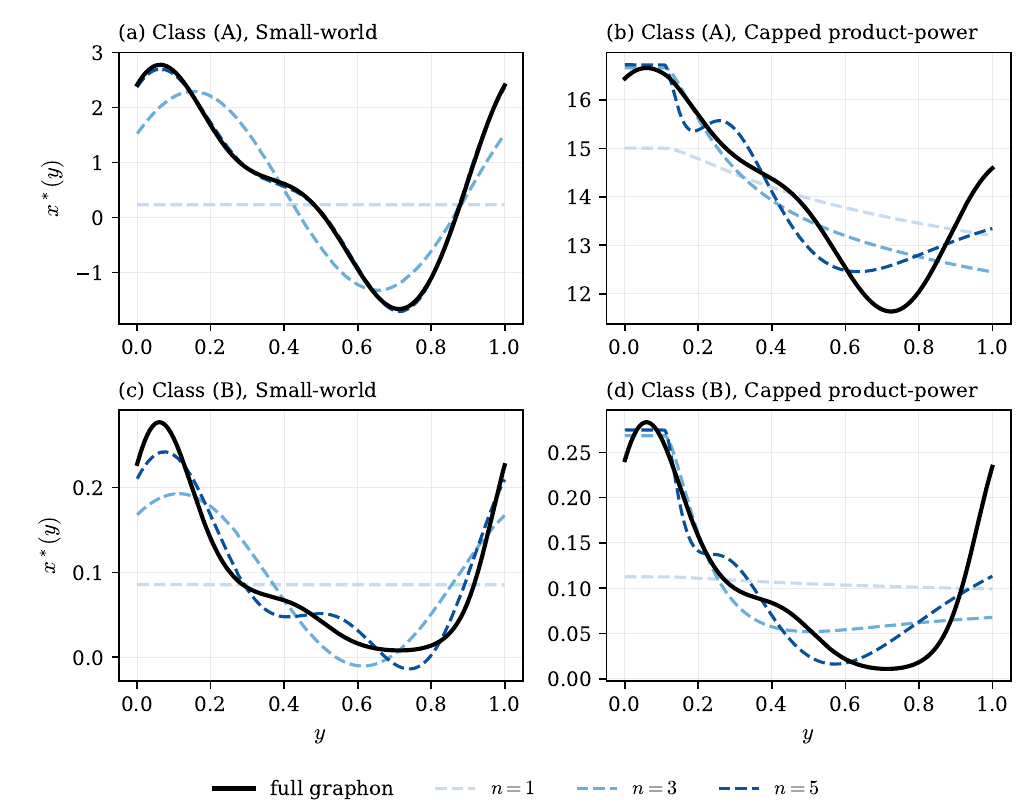}
\caption{Approximation of the equilibrium profile $x^\ast(y)$ for representative coupling strengths. (a) Class-\textup{(A)} dynamics on the small-world graphon, $D=15$. (b) Class-\textup{(A)} dynamics on the capped product-power graphon, $D=15$. (c) Class-\textup{(B)} dynamics on the small-world graphon, $D=8$. (d) Class-\textup{(B)} dynamics on the capped product-power graphon, $D=4$. Each panel compares the full equilibrium profile with spectral reconstructions using $n=1,3$, and $5$ variables.}
\label{fig:profiles-num}
\end{figure}

Figure~\ref{fig:profiles-num} shows the same approximation in the state space. It compares the full equilibrium profile $x^\ast(y)$ with its $n$-eigenfunction reconstructions for $n \in \{ 1, 3, 5\}$. For each pair of graphon and dynamical system, we set the coupling strength, $D$, such that the equilibrium is spatially heterogeneous: $D=15$ for the class-\textup{(A)} dynamics on the small-world and capped product-power graphons (Figs.~\ref{fig:profiles-num}(a) and (b)), $D=8$ for the class-\textup{(B)} dynamics on the small-world graphon (Fig.~\ref{fig:profiles-num}(c)), and $D=4$ for the class-\textup{(B)} dynamics on the capped product-power graphon (Fig.~\ref{fig:profiles-num}(d)). In all four cases, we visually observe that the approximation error becomes smaller as $n$ increases.
In particular, in Fig.~\ref{fig:profiles-num}(a), the equilibrium of the reduced dynamical system with $n=5$ is almost indistinguishable from the full profile.

\section{Conclusions}\label{sec:conclusions}

We developed the spectral dimension reduction of graphon dynamical systems as an operator-theoretic counterpart of finite-network spectral reduction. The central point is that, owing to its symmetry, a graphon defines a compact self-adjoint integral operator, so its eigenfunctions provide the natural reduction basis and its singular values are the absolute values of its eigenvalues. Using this structure, we formulated the $n$-dimensional spectral reduction directly at the graphon level and proved two types of results. First, for finite-rank graphons and compatible vector fields, the reduction is exact: the observable dynamics closes in dimension equal to the graphon rank (Theorem~\ref{thm:exact}), linear graphon dynamics decouples exactly into spectral modes (Corollary~\ref{cor:linear-exact}), and equilibria, linear stability, and local bifurcations are preserved by the reduced system (Theorem~\ref{thm:bifurcation}). Second, for infinite-rank graphons, the reduction is approximate but controlled. For class \textup{(B)} graphon dynamical systems, the discarded interaction modes enter the error estimate through the spectral tail $\abs{\lambda_{n+1}}$ (Theorem~\ref{thm:align}). For more general graphon dynamical systems, the error is controlled by the state truncation error (Proposition~\ref{prop:align-general}); the resulting observable trajectory error is bounded by a Gr\"onwall estimate (Theorem~\ref{thm:accuracy}). These results make precise, in the graphon setting, how low rank or rapid spectral decay enables low-dimensional descriptions of network dynamics.

This paper builds on our previous study of graphon dynamical systems~\cite{Eldo2025arxiv}, which established the continuum limit of a class of nonlinear dynamical systems on networks and proved existence, uniqueness, and boundedness of solutions of the limiting graphon dynamical system given by~\eqref{eq:graphon-ds}. That work then focused on the graphon version of the Gao--Barzel--Barab\'asi (GBB) reduction~\cite{gao2016universal}, a one-dimensional reduction method for dynamical systems on networks, and proved that the finite-network GBB observable converges to its graphon counterpart. Such convergence of observables, however, does not by itself imply that the GBB-reduced system accurately represents the full dynamics. Indeed, except in special cases, it is generally too much to expect a single reduced variable to provide a quantitatively accurate approximation of a high-dimensional dynamical system~\cite{gao2016universal, kundu2022accuracy}, although the GBB reduction often successfully captures bifurcation and resilience structures~\cite{gao2016universal, tu2017collapse, kundu2022accuracy, masuda2022dimension}. In the present paper, we instead focused on spectral reduction, which is systematically based on spectral/SVD structure and allows multiple reduced variables. We established exactness and error-control results for the spectral reduction of graphon dynamical systems at arbitrary reduced dimension $n$.

Several further problems remain.
First, Theorem~\ref{thm:bifurcation} concerns local bifurcations of equilibria; the corresponding statement for global bifurcations, and for dynamics that are only approximately of class \textup{(A)} or \textup{(B)}, is not addressed here. Second, relating the spectral reduction of a large finite network to the graphon reduction studied here would require the convergence of the eigenvectors of the network matrices to the eigenfunctions of $\Tw$, which is known only under spectral-gap and regularity conditions~\cite{gong2022regular, borgs2012convergent}. We regard establishing it for general graphons as the main open problem adjacent to this work. Third, extensions of the theory to sparse graphons~\cite{kaliuzhnyi2017semilinear}, to directed (non-symmetric) graphons where the genuine SVD with distinct left and right singular functions is needed, and to multidimensional node states, are natural directions for future research.

\section*{Acknowledgments}

N.M.\ acknowledges financial support by the Japan Science and Technology Agency (JST) Moonshot R\&D (under Grant No. JPMJMS2021), the National Science Foundation (under grant no.\,2204936), and JSPS KAKENHI (under grant nos.\,JP 23H03414, 24K14840, and 24K03013).
During the preparation of this manuscript, the authors used Claude Opus 4.7, University of Michigan (U-M) ChatGPT 5.5, and ChagGPT 5.5 Pro for language editing, code development, and mathematical analysis assistance. The authors manually verified all code and manuscript text.

\bibliographystyle{unsrtinit}
\bibliography{references}

\end{document}